\documentclass[a4paper,11pt,dvipdfm]{article} 

\usepackage{amsmath} 
\usepackage{amssymb} 
\usepackage{amsthm} 

\usepackage{epsfig} 
\usepackage{caption2} 
\usepackage[import,frame]{xy} 
\input ./warmread.sty 
\let\xyWARMprocess\xyWARMprocessMo 
\let\WARMprocessEPS\WARMprocessMoEPS 
\wrmquiet 
\usepackage[pagebackref,colorlinks=false]{hyperref} 
\newcommand\MR[1]{\relax\ifhmode\unskip\spacefactor3000 \space\fi
\MRhref{#1}{MR #1}} 

\newcommand{\etal}{\emph{et~al.}} 

\newtheorem{theorem}{Theorem}[section] 
\newtheorem{lemma}[theorem]{Lemma} 
 
\newtheorem{corollary}[theorem]{Corollary} 

\theoremstyle{definition} 
\newtheorem{example}[theorem]{Example} 
 
\newtheorem{definition}[theorem]{Definition}

\newcommand{\id}{\mathrm{id}} 
\newcommand{\real}{\mathbb{R}} 
\newcommand{\nat}{\mathbb{N}} 
\newcommand{\integer}{\mathbb{Z}} 
\newcommand{\rational}{\mathbb{Q}} 
\newcommand{\rat}{\mathbb{Q}} 
 
\newcommand{\lattice}{\mathcal{L}}

\newcommand{\ol}{\overline} 
 
\newcommand{\tile}{\ensuremath{\mathcal{T}}}

\setcaptionmargin{18pt}

\begin{document}
	
	\author{Edmund O. Harriss$^1$ and Jeroen S. W. Lamb$^2$
\\\\
 \small{$^1$ School of Mathematical Sciences, Queen Mary University of London,
 London E1 4NS, UK} \\
\small{edmund.harriss@mathematicians.org.uk}
\\
 \small{$^2$ Department of Mathematics, Imperial College London,
 London SW7 2AZ, UK}
 }
	\title{One-dimensional substitution tilings with an interval projection structure} 
	\date{\today} 
	\maketitle 
	\begin{abstract}
	We study nonperiodic tilings of the line obtained by a projection method with an \emph{interval projection structure}.  We obtain a geometric characterisation of all interval projection tilings that admit substitution rules and describe the set of substitution rules for each such a tiling. We show that each substitution tiling admits a countably infinite number of nonequivalent substitution rules.  We also provide a complete description of all tilings of the line and half line with an interval projection structure that are fixed by a substitution rule.   Finally, we discuss how our results relate to renormalization properties of interval exchange transformations (with two or three intervals).  
	\end{abstract}
	\tableofcontents\normalsize 

\newpage	
	\section{Introduction} 
	
	The study of substitution rules for symbolic sequences and one-dimensional tilings has a long and rich history. In particular, Sturmian sequences and their substitution rules have been intensively studied in the context of computer science and dynamical systems \cite{Fogg:SIDAC, Lothaire:ACOW}.    
	
	Hedlund and Morse \cite{Morse:SD2ST} showed that Sturmian sequences can be given an elegant geometric construction (see also \cite{Series:GMCOG,Series:TGOMN}): every Sturmian sequence can be obtained from the intersections of a straight line with a square grid by assigning the symbol 0 to each intersection of the line with a horizontal grid line and a 1 to each intersection with a vertical grid line. The binary {\em cutting} sequence thus obtained is a Sturmian sequence. Sturmian substitution rules are those morphisms of binary strings that map Sturmian sequences to Sturmian sequences. These substitution rules have an elegant geometric realization: a unimodular matrix with integer coefficients inducing a linear transformation of the line associated with the initial Sturmian sequence to the line associated with its image under the substitution rule. Interestingly, the Sturmian substitution rules are generated by only three morphisms.  Incidentally these three morphisms also generate the automorphism group of the free group on two symbols $F_2$, see eg \cite{Peyriere:EDCAI,Lamb:OTCPM,Wen:LIOIS}. 
	
	In the physics literature, in the area of quasicrystals \cite{Senechal:QAG,Axel:BQ}, Sturmian sequences arise as symbolic representations of \emph{canonical projection tilings}. Their construction, which we describe in Section \ref{sect:Main_results}, is closely related to the one by Hedlund and Morse mentioned above. 

	In this paper we consider one-dimensional \emph{interval projection tilings}. This class of tilings appears to have been first considered in \cite{Masakova:SRFAS}, and contains the above mentioned set of canonical projection tilings. 
	
	Our results unify and extend various results in the literature. In particular, we answer the following central questions: 
	\begin{itemize}
		\item Which one-dimensional tilings with an interval projection structure are substitution tilings? 
		\item What substitution rules do such tilings admit? 
	\end{itemize}
	We note that the first problem was previously only resolved for one-dimensional canonical projection tilings. The second question had not been addressed even in the case of particular well-studied examples such as the \emph{Fibonacci} tiling.
	
	Whereas many studies of nonperiodic tilings have been mainly combinatorial or algebraic, in this paper we take a predominantly geometric point of view.
	
	Importantly, it turns out that the results of this paper provide some key insights enabling the understanding of the existence of substition rules for higher dimensional (canonical) projection tilings, including the well-known Penrose and Ammann tilings \cite{Harriss:OCST,Harriss:CSTOA} which serve as a prototypical examples of quasicrystalline tilings. 

	\section{Main results} \label{sect:Main_results} 
	
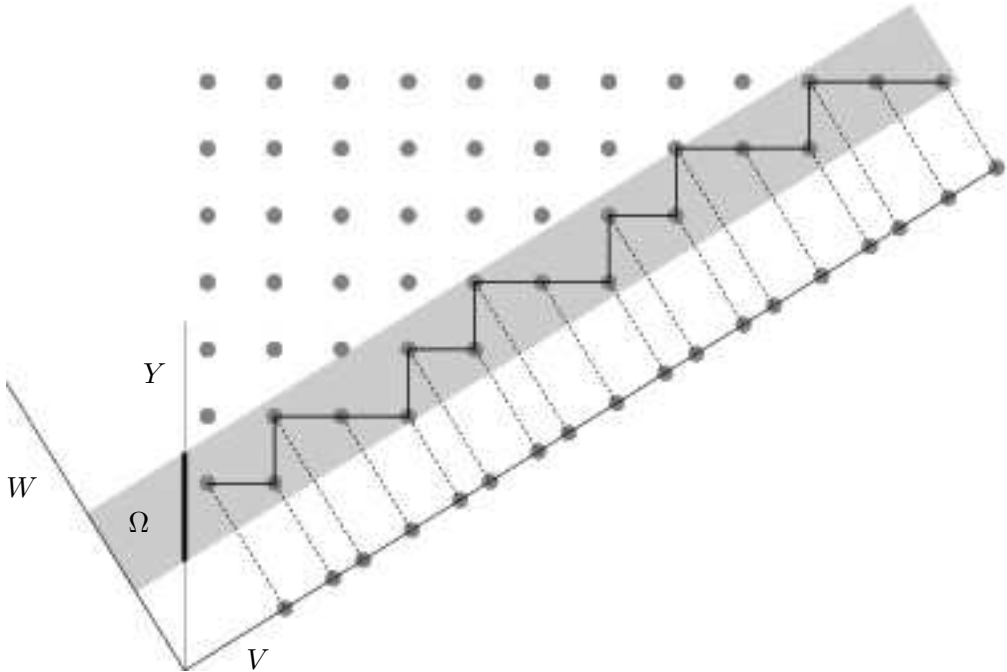
\begin{figure}[htp]
	$$
	\begin{xy}
		\xyMarkedImport{}
		\xyMarkedPos{1}*\txt\labeltextstyle{$V$}
		\xyMarkedPos{2}*\txt\labeltextstyle{$Y$}
		\xyMarkedPos{3}*\txt\labeltextstyle{$W$}
		\xyMarkedPos{4}*\txt\labeltextstyle{$\Omega$}
	\end{xy}
	$$
	\caption[Fibonacci tiling of the line]{Illustration of the geometric structure of the Fibonacci tiling, which is the prototypical example of a one-dimensional tiling with an interval projection structure. Consider a strip in $\real^2$, indicated in grey and defined by the translation of a unit square in the direction of a line $V$ with slope $\tau^{-1}$, where $\tau = (1+\sqrt{5})/2$ denotes the golden ratio. The intersection of $\integer^2$ with this strip is projected down the space $W$ to $V$ to yield the set of vertices of the Fibonacci tiling on $V$.  The window is the interval of the space $Y$ that generates the strip when translated along $V$.} \label{fig:construction}
	\end{figure}

	 We consider tilings of the line by intervals obtained by the following geometrical construction. Let $V$ and $W$ be two transversally intersecting lines in $\real^2$. We consider the intersection of the lattice $\integer^2\subset\real^2$ with a strip $V+\Omega\subset\real^2$, where $\Omega$ is a half open interval of the vertical axis, referred to as the {\em window}. For later convenience, we adopt the convention that the window lies in the subspace $Y$ parallel to the direction the vertical lattice generator. Subsequently, we project the lattice points inside $\integer^2\cap(V+\Omega)$ to $V$, choosing the projection $\Pi_V$ \emph{parallel to $W$} (such that $\Pi_V^{-1}(v)$ is a line parallel to $W$ for all $v\in V$). The resulting discrete point set $$ \Pi_V(\integer^2\cap(V+\Omega))\subset V $$ is considered as the set of vertices for a tiling of the line $V$, where tiles are intervals. The construction is sketched in Figure~\ref{fig:construction}. We refer to the resulting tilings as {\em one-dimensional tilings with an interval projection structure}. In this construction we always assume that $V$ and $W$ are not parallel to one of the coordinate axes (thus avoiding some degenerate constructions yielding periodic tilings). 
	
	\WARMprocessEPS{introduction_can2to1}{eps}{bb}
	
	In the special case where the length of the window $\Omega$ is precisely equal to the length of the projection $\Pi_Y$ to $Y$ (parallel to $V$) of a $1\times 1$ unit square in $\real^2$, the points inside the strip can be connected by horizontal and vertical line pieces ({\em steps}) to form a monotone {\em staircase}. This staircase can be viewed as a {\em discrete} approximation of a line parallel to $V$. Such staircases can be represented by sequences of two symbols, each symbol representing a type of step (horizontal or vertical). Such binary sequences are the {\em Sturmian} sequences and have been extensively studied in the literature \cite{Lothaire:ACOW}, also in relation to substitution rules \cite{deBruijn:SOZAO,Seebold:OTCOS,Ito:OCFSA,Parvaix:PDDMS,Berstel:RRISW,Crisp:SICS}.  Sturmian sequences arise in number theory as the sequence of differences between consecutive terms in a Beatty sequences (the numbers $(\lfloor(n \lambda\rfloor))_{n=1}^\infty$ where $\lambda$ is irrational).  For discussions of substitution rules in this context see \cite{Stolarsky:BSCFA,Fraenkel:DONTB,deBruijn:UGOBS,Komatsu:SIBS,Komatsu:SIIBS}.
		
	When $W$ is chosen such that the projections parallel to $W$ onto $V$ of consecutive steps of the staircase intersect only at their endpoints, the tiling of $V$ can be represented by a binary sequence; each symbol representing projection of either a horizontal or a vertical step. 
	
	In general, we consider tilings of the line by intervals (of finite length), where the intervals are taken from a finite {\em protoset} of tiles (sometimes called {\em prototiles}). Every tile in the protoset has a geometric shape (an interval with a given length) and possibly a {\em label}. A tiling consists of an infinite collection of translated copies of tiles from the protoset that cover the line in such a way that translated prototiles only overlap on their boundaries. Labels provide a way to distinguish between geometrically identical tiles in a protoset. 
	
	In this paper we provide a characterisation of all one-dimensional tilings with an interval projection structure that admit a substitution rule. We also characterise the set of all substitution rules that such a tiling admits. 
	
	The (geometric) substitution rules that we consider consist of two parts. First, we perform a uniform inflation of each (labeled) tile by a constant factor, say $\lambda>1$. Thereafter, we replace each inflated tile by a set of uninflated tiles in such a way that each inflated tile with the same label is replaced in exactly the same way. Here, we always require that the substitution is {\em vertex hierarchic}, implying in the present context that each inflated tile is exactly covered by the set of tiles that replace it. We say that a tiling is a substitution tiling if the image of a tiling under a substitution rule is {\em locally isomorphic} to the original tiling. Two tilings are locally isomorphic if they are on any local scale indistinguishable in the sense that any finite patch of the one also appears in the other. A set of locally isomorphic tilings is called a {\em local isomorphism class}. In the case of one-dimensional tilings with an interval projection structure, local isomorphism classes are characterized by the spaces, $V$ and $W$ and the window length $|\Omega|$ (see Lemma~\ref{2to1_LI_class}).
	
\begin{figure}
		[htp] $$ 
		\begin{xy}
			\xyMarkedImport{} 
			\xyMarkedPos{1}*\txt\labeltextstyle{(a)} 
			\xyMarkedPos{2}*\txt\labeltextstyle{(b)} 
			\xyMarkedPos{3}*\txt\labeltextstyle{(c)} 
			\xyMarkedPos{4}*\txt\labeltextstyle{(d)}  
		\end{xy}
		$$ \caption[The Fibonacci substitution rule window]{Geometric derivation of a substitution rule for the Fibonacci tiling.\\
		(a) A patch of the Fibonacci tiling and the corresponding staircase. The tiling space $V$ and $W=V^\perp$ (along which we project to $V$) are the expanding and contracting eigenspaces of the matrix 
$M=\left( 
	\begin{array}
		{cc}1&1\\1&0 
	\end{array}
	\right)$.  
\\
		(b) The vertices of the staircase that have a vertical step to their right are identified.  They correspond to lattice points lying in the dark grey strip.\\
		(c) After removal of these points, we obtain a staircase with horizontal and diagonal steps.  We may interpret this removal as the substitution 
$ab\to a, a\to b$, where $a$ represents a long and $b$ represents a short tile. We note that in the orginal tiling horizontal staircase steps project to long tiles and vertical steps to short tiles, whereas in the new staircase horizontal steps project to short tiles and diagonal steps to long tiles.\\
		(d) By application of the lattice automorphism $M^{-1}$, the new staircase is mapped to a staircase for a Fibonacci tiling (as characterized by the slope of $V$ and the same window length), as it commutes with the projections along $V$ and along $W$.  It follows that the tiling for the new staircase is a Fibonacci tiling that is uniformly inflated by a factor equal to the expanding eigenvalue of $M$.  The substitution rule  $a\to ab, b\to a$ is known as the Fibonacci substitution rule and maps Fibonacci tilings to Fibonacci tilings.}
	 \label{fig:fibsub} 
	\end{figure}
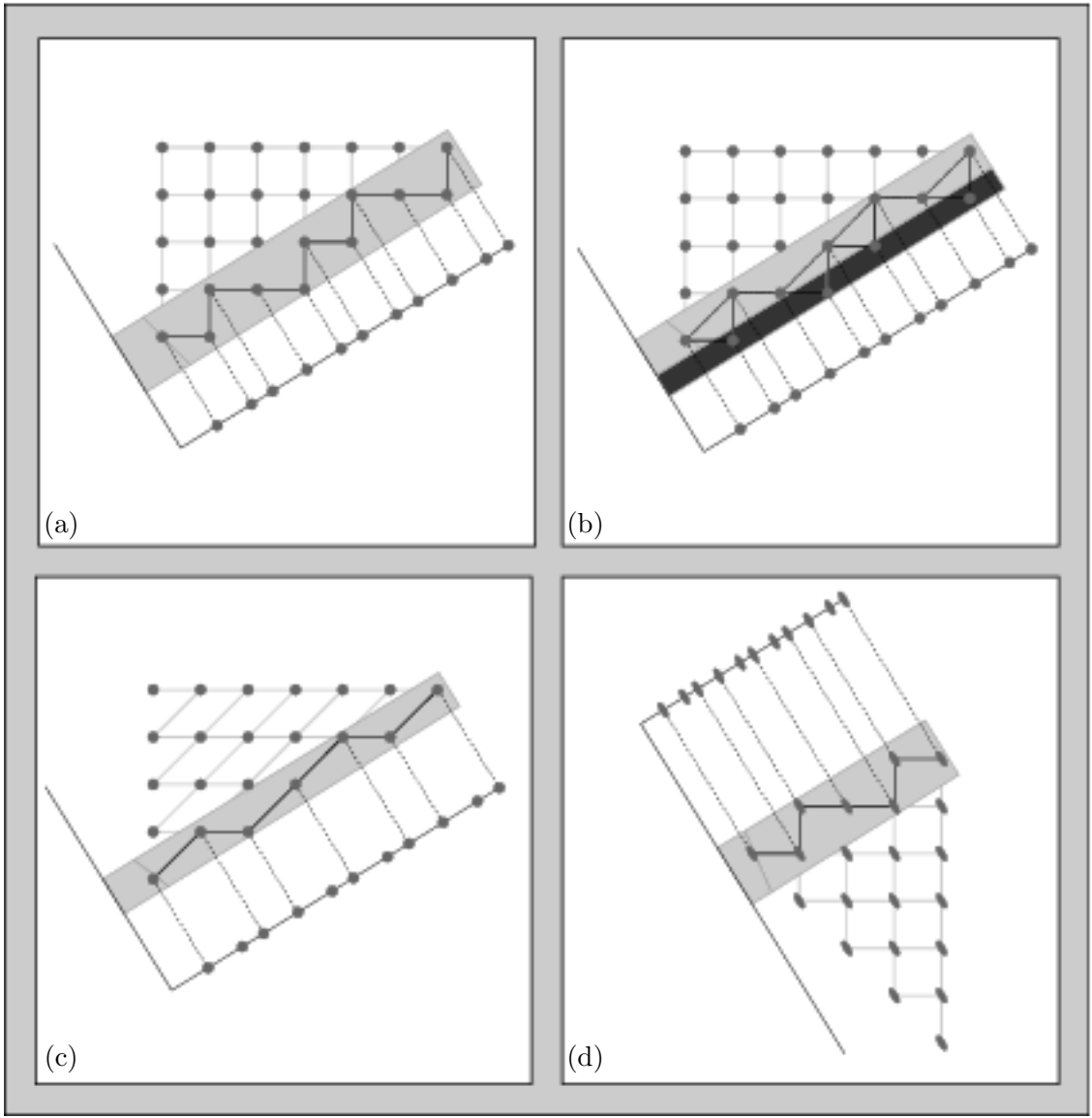

	A well known example of a one-dimensional substitution tiling with an interval projection structure is a Fibonacci tiling. Its construction is as sketched in Figure \ref{fig:construction} with $V$ having its slope equal to the inverse of the golden ratio ($\tau^{-1}$), $W=V^\perp$, and the window 
$\Omega$ having  length $1+\tau^{-1}$. The corresponding tiling can be represented by a binary Sturmian sequence and admits the well known Fibonacci substitution rule. A geometric derivation of the Fibonacci substitution rule is presented in Figure~\ref{fig:fibsub}.  

	\WARMprocessEPS{introduction_fibbwindow}{eps}{bb} 
	
	In this paper we establish how the geometric approach to a substitution rule for the Fibonacci tiling, as illustrated in Figure~\ref{fig:fibsub}, generalizes to provide a full characterization of all substitution rules for one-dimensional substitution tilings with an interval projection structure. 

It should be noted that we view the tilings in the first instance as geometric objects (a partition of the line by intervals) without considering labels.  Such a geometric tiling admits a substitution rule if there exists a labelling of the tiles (respecting the tile-shapes, ie two intervals of different length cannot be assigned the same label), such that the labelled tiling admits a substitution rule\footnote{In particular this rules out certain Sturmian sequences that are substitutive with respect to symbolic substitution rules (ie they can be generated by projection on the symbols).  This is because the projection breaks the geometric requirement that tiles of different lengths cannot have the same label.  For example consider the substitution rule $a\to ac, b\to acacb, c\to acacb$.  When this substitution rule is brought into the geometric setting the tiles come in two lengths, the $a$ tiles and the $b$ and $c$ tiles.  As the substitution rule does not differentiate between $b$ and $c$, we can remove the labels to obtain a Sturmian sequence that can be constructed using a geometric substitution rule.  However if we bring together the $a$'s and the $b$'s we get a different Sturmian sequence.  This symbolic projection, however, gives the same label to tiles of different lengths, thus we do not consider this to be a geometric substitution rule.  For a general discussion of substitutive Sturmians from a symbolic point of view see~\cite{Arnoux:RGSC2}.}.
	
	The first result of this paper is summarized in the following theorem: 
	\begin{theorem}
		\label{thm:char} Let $V$ and $W$ be the tiling and window spaces for a nonperiodic one-dimensional tiling with an interval projection structure. This projection tiling admits a substitution rule if and only if $V$ and $W$ are the expanding and contracting eigenspaces of a primitive matrix in $Gl(2,\integer)$, and the window length is inside $\rational[\lambda]$ where $\lambda$ is the expanding eigenvalue of $M$.\footnote{Recall that $\rational[\lambda]=\rational+\lambda\rational$ is the smallest algebraic field containing $\lambda$, $Gl(2,\integer)$ is the group of $2\times 2$ matrices with integer coefficients whose determinant is equal to $\pm 1$.  A matrix $M \in Gl(2,\integer)$ is \emph{primitive} if their exists $n \in \nat$ such that all entries of $M^n$ are non-zero.} 
	\end{theorem}
	It is important to note that the above theorem deals with {\em nonperiodic} tilings. Periodic tilings are special, as these may admit substitution rules that bare no relationship with the geometric setting of Theorem~\ref{thm:char}. For instance, if a one-dimensional tiling with an interval projection structure is periodic with unit cell $U$, and $W$ is chosen such that the ratios of the lengths of the tiles are rational, then there exists an appropriate inflation (multiplication by some integer) such that all inflated tiles can be replaced by a number of concatenated unit cells $U$. 
	
	In Theorem~\ref{thm:char} no reference is made to the protoset. In fact, it turns out that each one-dimensional substitution tiling with an interval projection structure admits a countable infinity of substitution rules, using many different labelings of the tiles (and thus different protosets). The following theorem provides a geometrical characterization of all substitution rules of one-dimensional substitution tilings with an interval projection structure. 
	\begin{theorem}
		\label{thm:subchar} Consider a nonperiodic one-dimensional substitution tiling with interval projection structure and window $\Omega$, as in Theorem~\ref{thm:char}. Then for each choice of the (half open) interval window $\Omega'\subset W$ such that 
		\begin{itemize}
			\item[(i)] $|\Omega'|=|\Omega|$, 
			\item[(ii)] $M(\Omega+V)\subset (\Omega'+V)$,
			\item[(iii)] the length of each connected component of 
$\Omega'\setminus \Pi_Y M(\Omega)$ is in $\rational[\lambda]$, 
		\end{itemize}
		there exists a substitution rule that maps the projection tiling with window $\Omega$ to the locally isomorphic projection tiling with window $\Omega'$. 
		
		Moreover, for each substitution rule one can identify a window $\Omega'$ with the above mentioned properties. 
	\end{theorem}
	
	Among the substitutions mentioned in Theorem~\ref{thm:subchar}, we may distinguish between \emph{local} and \emph{nonlocal} substitutions. Local substitutions are those where the shape of a finite neighbourhood of a tile (a patch consisting of a finite number of adjacent tiles) determines its label. In many studies of substitution rules, such as \cite{Pleasants:DQ,Baake:QTWTS}, the locality of substitution rules is assumed. However, the following result shows that although one-dimensional substitution tilings with an interval projection structure always admit infinitely many nonlocal substitution rules, only some admit local substitutions.

	\begin{theorem}\label{thm:local} 
Consider the characterization of substitution rules in Theorem~\ref{thm:char}
and Theorem~\ref{thm:subchar}. Let the window pair $\Omega$ and $\Omega'$
be associated to a substitution rule for a tiling with interval projection structure. Let $M\in Gl(2,\integer)$ be the primitive matrix 
associated to this substitution rule. Moreover, let $d_+$ and $d_-$ denote the 
distances between the tops and bottoms, respectively, of the windows 
$\Pi_Y M(\Omega)$ and $\Omega'$.  

Then the substitution rule is local if and only if 
at least one of the following conditions is satisfied:
\begin{itemize}
\item $(\det(M)-\lambda)|\Omega|, d_+,d_-\in \integer[\lambda]$,
\item $(\det(M)-\lambda)|\Omega|, \lambda d_+ - \det(M)|\Omega|, 
        \lambda d_- - \det(M)|\Omega| \in \integer[\lambda]$.
\end{itemize}
\end{theorem} 
\begin{corollary}
A substitution tiling with interval projection structure whose window length
satisfies
$$
	|\Omega| \notin \frac{\integer[\lambda]}{\det M+\lambda}\cup \frac{\integer[\lambda]}{\det M-\lambda})
$$
only admits nonlocal substitution rules. \label{cor:nonloc}
\end{corollary}
	To the best of our knowledge, the existence of nonlocal substitution rules was not known before, not even for the widely studied Fibonacci tilings, for which we provide an example. 

\begin{example}
		[Nonlocal substitution rule for Fibonacci tilings]\label{ex:nonlocalfib} Let the tiling space $V$ and projection direction $W$ be the expanding and contracting eigenspaces of the matrix 
\begin{equation}	M=\begin{pmatrix}
			1&1\\
			1&0 
		\end{pmatrix}.\label{fibmatrix}\end{equation}
		The Fibonacci tiling is obtained by the projection method, with a window that is obtained as the projection  of a unit square along $V$ to $Y$. 
		
		We consider the substitution rule where the expansion predecessor is associated with the subwindow $M(\Omega)$ that lies precisely at the centre of the initial window $\Omega'$. It is readily verified that the distance between the endpoints of the subwindow $M(\Omega)$ and the endpoints of $\Omega'$ is $\tau^{-1}/2 \in\rat[\lambda]\setminus(\integer[\lambda])$. Hence, by Theorem~\ref{thm:local} the corresponding substitution rule is nonlocal. It is defined on a protoset consisting of four tiles: two long intervals ($a_1$ and $a_2$) and two short intervals ($b_1$ and $b_2$). It terms of these, the substitution rule takes the form 
		$$\left\{ 
		\begin{array}
			{l}a_1\to b_1 a_1,\\a_2\to a_2 b_2,\\b_1\to a_2,\\b_2\to a_1. 
		\end{array}
		\right.$$ The derivation of this substitution rule is presented in Example~\ref{ex:nonlocalfibdet}.
\end{example}
\begin{example}[Substitution tiling with only nonlocal substitution rules]
Consider a projection tiling with $V$ and $W$ the expanding and contracting 
eigenspaces of the matrix $M$ of (\ref{fibmatrix}), as in the case of the Fibonacci tiling. We consider tilings with window length $1/n$ with $n>1$ integer.
We have $\lambda=-\tau$ and $\det M=-1$ and $1/n\not\in \integer[\tau]/(1+\tau)\cup \integer[\tau]/(1-\tau)$,
so that by Theorem~\ref{thm:char} and
Corollary~\ref{cor:nonloc} the tiling admits substitution rules but none of them is a local substitution rule.
\end{example} 
	Historically, there has been a considerable interest in tilings that are \emph{fixed} by a substitution in the sense that the tiling, after the application of a substitution rule, is an exact (translated) copy of the original tiling.\footnote{We adopt the convention that two tilings are equivalent if they can be mapped onto each other by a translation.\label{fn:equivalence}}
	
	The following theorem summarizes our results about the existence of one-dimensional substitution tilings with an interval projection structure that are fixed by a substitution rule. 
	\begin{theorem}
		\label{introthm:fix_tilings} Consider the set of substitution rules $\mathcal{S}$ fixing the local isomorphism class of a one-dimensional substitution tiling with an interval projection structure associated to a window length $|\Omega|$ and matrix $M\in Gl(2,\integer)$. Then, a substitution rule in $\mathcal{S}$ fixes $\det(M-I)$, $\det(M-I) - 1$, or $\det(M-I) + 1$ tilings, where $I$ denotes the identity matrix. 
	\end{theorem}
 
	The conditions on a substitution rule that determine which one of the three options given in the above theorem arises, can be expressed in terms of properties of the windows $\Omega$, $\Omega'$ and the matrix $M$. We refer the reader to Theorem~\ref{thm:fixed}. In Section~\ref{sec:fixed} we also describe the set of one-dimensional tilings with an interval projection structure that tile the half-line $\real^+$ and are left invariant by a substitution rule. 
	
Our geometric methods are inspired by and adapted from the algebraic methods used by Masakova~\etal\ in \cite{Masakova:SRFAS}. Masakova~\etal\ present a large class of quasicrystal tilings of the line. Their tilings have an interval projection structure where $V$ and $W$ are the expanding and contracting eigenspaces of a primitive matrix $M$ with integer coefficients. Masakova~\etal\ use a different, but equivalent, description of these tilings in terms of algebraic number theory. Under the assumption that $M$ is unimodular they show that such a tiling admits a labelling and a substitution rule on the labelled protoset that fixes the tiling, if and only if both endpoints of the window lie in $\rational[\lambda]$ in $Y$. This can be viewed as a corollary of the results in this paper. 

Our results can also be viewed in the context of renormalization properties of two and three-interval exchange transformations. For a discussion, see Section~\ref{sec:iet}.

As the proofs of our results are all constructive, they can be viewed as the basis of an algorithm for constructing substitution rules. In Table~\ref{table_examples} we have listed some examples of substitution rules found using the methods described in this paper, encoded as a Mathematica program \cite{Harriss:NSTOT_code}. 

It turns out that the substitution rules in Table~\ref{table_examples} (and all others we obtained) are invertible, in the sense that it has a formal inverse so that it corresponds to an automorphism of a free group. It is well known that binary Sturmian substitution rules correspond to automorphisms of the free group $F_2$ \cite{Wen:LIOIS, Peyriere:EDCAI,Lamb:OTCPM, Ei:DTOIS}. Our experiments thus suggest that all substitution rules for one-dimensional tilings with an interval projection structure are invertible.
	
	\begin{table}
		[htp] \centering 
		\begin{tabular}
			{|c|c|ccc|c|c|} \hline
			
			$M$ & $l$ & $s_1$&$-$&$s_2$ & Substitution & Inverse  \\
			
			\hline \hline 
			
			$ 
			\begin{pmatrix}
				1&1\\2&1 
			\end{pmatrix}
			$& 
			
			$ 
			\begin{pmatrix}
				-1\\1 
			\end{pmatrix}
			$&
			
			$ 
			\begin{pmatrix}
				-\frac{1}{2}\\0 
			\end{pmatrix}
			$
			
			&$-$&$ 
			\begin{pmatrix}
				-\frac{1}{2}\\1 
			\end{pmatrix}
			$&$ 
			\begin{matrix}
				a\to b_1ab_2\\b_1\to b_1a\\b_2\to ab_2 
			\end{matrix}
			$&
			
			$ 
			\begin{matrix}
				a\to b_2a^{-1}b_1\\b_1\to ab_2^{-1}\\
				
				b_2\to b_1^{-1}a 
			\end{matrix}
			$ \\
			
			\hline
			$ 
			\begin{pmatrix}
				1&1\\2&1 
			\end{pmatrix}
			$& 
			
			$ 
			\begin{pmatrix}
				-1\\1 
			\end{pmatrix}
			$&
			
			$ 
			\begin{pmatrix}
				-\frac{1}{4}\\0 
			\end{pmatrix}
			$
			
			&$-$&$ 
			\begin{pmatrix}
				-\frac{1}{4}\\1 
			\end{pmatrix}
			$& 
			$ 
			\begin{matrix}
				
				a_1\to b_1b_1a_1\\
				
				a_2\to b_1a_3b_2\\
				
				a_3\to b_1a_3b_1\\
				
				b_1\to b_1a_2\\
				
				b_2\to b_1a_1 
			\end{matrix}
			$&
			
			$ 
			\begin{matrix}
				
				a_1\to b_2a_1^{-1}b_2\\
				
				a_2\to b_2a_1^{-1}b_1\\
				
				a_3\to b_2a_1^{-1}a_3b_2a_1^{-1}\\
				
				b_1\to a_1b_2^{-1}\\
				
				b_2\to a_1b_2^{-1}a_3^{-1}a_2 
			\end{matrix}
			$\\
			
			\hline
			$ 
			\begin{pmatrix}
				1&1\\2&1 
			\end{pmatrix}
			$& 
			
			$ 
			\begin{pmatrix}
				-\frac{1}{2}\\1 
			\end{pmatrix}
			$&
			
			$ 
			\begin{pmatrix}
				0\\0 
			\end{pmatrix}
			$
			
			&$-$&$ 
			\begin{pmatrix}
				-\frac{1}{2}\\0 
			\end{pmatrix}
			$& 
			$ 
			\begin{matrix}
				
				a\to bc\\
				
				b\to ba\\
				
				c\to bcc 
			\end{matrix}
			$&
			
			$ 
			\begin{matrix}
				
				a\to a^{-1}ca^{-1}b\\
				
				b\to ac^{-1}a\\
				
				c\to a^{-1}c 
			\end{matrix}
			$\\
			
			\hline
			$ 
			\begin{pmatrix}
				1&1\\2&1 
			\end{pmatrix}
			$& 
			
			$ 
			\begin{pmatrix}
				-\frac{1}{2}\\1 
			\end{pmatrix}
			$&
			
			$ 
			\begin{pmatrix}
				-\frac{1}{2}\\0 
			\end{pmatrix}
			$
			
			&$-$&$ 
			\begin{pmatrix}
				-\frac{1}{2}\\
				\frac{1}{2} 
			\end{pmatrix}
			$ &
			$ 
			\begin{matrix}
				
				a\to cb\\
				
				b\to c\\
				
				c\to cbab 
			\end{matrix}
			$&
			
			$ 
			\begin{matrix}
				
				a\to a^{-1}ca^{-1}b\\
				
				b\to b^{-1}a\\
				
				c\to b 
			\end{matrix}
			$\\
			
			\hline
			$ 
			\begin{pmatrix}
				2&3\\1&1 
			\end{pmatrix}
			$& 
			
			$ 
			\begin{pmatrix}
				-1\\\frac{2}{3} 
			\end{pmatrix}
			$&
			
			$ 
			\begin{pmatrix}
				0\\0 
			\end{pmatrix}
			$
			
			&$-$&$ 
			\begin{pmatrix}
				0\\
				\frac{1}{3} 
			\end{pmatrix}
			$ &
			$ 
			\begin{matrix}
				
				a_1\to ba_1a_2\\
				
				a_2\to ca_2\\
				
				b\to ca_1a_2\\
				
				c\to ca_1ca_1a_2 
			\end{matrix}
			$&
			
			$ 
			\begin{matrix}
				
				a_1\to bc^{-1}ba_2^{-1}cb^{-1}\\
				
				a_2\to bc^{-1}b\\
				
				b\to a_1b^{-1}a_2b^{-1}cb^{-1}\\
				
				c\to a_2b^{-1}cb^{-1} 
			\end{matrix}
			$\\
			
			\hline 
		\end{tabular}
		
		\caption{Some examples of substitution rules for one-dimensional tilings with an interval projection structure, where $V$ and $W$ are the expanding and contracting eigenspaces of a matrix $M\in Gl(2,\integer)$, and 
the width of the window $\Omega$ is represented by $l \in \rational^2$: $|\Omega| = \Pi_Y(l)$. The subwindow associated with the substitution rule is represented by the vectors $s_1,s_2\in\rational^2$: $\Pi_Y(s_1)$ and $\Pi_Y(s_2)$ are the relative positions of the bottom and top of the subwindow. Note that $M l=s_1-s_2$. The labels for the tiles are chosen such that those labelled $a_i$ corresponds to the projection of a horizontal staircase step, those labelled $b_i$ correspond to the the projection of a vertical staircase step and those labelled $c_i$ correspond to the projection of a diagonal staircase step (the combination of a horizontal and vertical step). In the last column we provide the formal inverse of the substitution rule.}
		\label{table_examples} 
	\end{table} 
	
	\section{Preliminaries} \label{background} 
	
	In this section we recall some definitions and properties of substitution rules and interval projection tilings. 
	
	\subsection{Substitution rules} \label{sec:find_subst} 
	
	Consider a set of prototiles $\tau=\{\tau_i\}$, and a real number $\lambda > 1$. Let $\lambda \tau= \{\lambda \tau_i\}$ denote the set of tiles obtained from $\tau$ after inflating each prototile by a uniform factor $\lambda$. A (vertex hierarchic) \emph{replacement rule} $\sigma_R$ associates a patch of tiles to each expanded prototile in $\lambda \tau$, so that $\sigma_R(\lambda \tau_i)$ covers the inflated prototile $\lambda \tau_i$ exactly. Consequently, the replacemt rule is \emph{vertex hierarchic}: the vertices satisfy the inclusion $vert(\tau_i)\subset vert \sigma_R(\lambda \tau_i)$.
	
	A \emph{substitution rule} $\sigma$ consists of a pair $(\lambda, \sigma_R)$, indicating the subsequent application of an uniform inflation by the factor $\lambda$ and a replacement rule $\sigma_R$. The action of a substitution rule $\sigma$ on a protoset $\tau$ induces, a substitution rule on any tiling (or patch) $\tile$ admitted by this protoset.  A patch or tiling is substituted by applying the uniform inflation, giving $\lambda\tile$, and replacing all translated copies of the inflated prototiles according to the replacement rule. 
	
	A tiling $\tile'$ is called a \emph{predecessor (tiling)} for $\tile$ if $\sigma(\tile') = \tile$. We refer to the inflated tiling $\lambda \tile'$ as the \emph{expansion predecessor (tiling)}. 
	
	\begin{definition}
		[Substitution Tiling] A tiling $\tile$, with protoset $\tau$, is a substitution tiling with substitution rule $\sigma$ if it has infinite predecessors under $\sigma$, ie if for all $n\in\nat$ there exists a locally isomorphic tiling $\tile'$ such that $\sigma^{n}\tile'=\tile$. 
	\end{definition}
	It follows that the set of substitution tilings for a primitive\footnote{A substitution rule is called \emph{primitive} if after a finite number of subsequent applications each tile (in the protoset) is replaced by a set of tiles including copies of all tiles in the protoset. As a consequence, under the subsequent application of a primitive substitution rule no patches grow up that contain only a subset of tiles from the protoset.} substitution rule forms a set of local isomorphism classes.\footnote{Note that in order to render this local isomorphism class unique, other (less intuitive) definitions of substitution tilings are used in the literature, see for example \cite{Goodman:MRAST}.  It turns out that the projection structure used in this paper forces the substitution rules to be associated to a unique local isomorphism class.}
	
	The geometric definition of substitution tiling given here is motivated by the study of geometric substitution rules in higher dimensions, see \cite{Senechal:QAG,Goodman:MRAST}.  In one dimension these are closely related to symbolic substitution rules and morphisms as studied in computer science and dynamical systems, see \cite{Allouche:AS,Fogg:SIDAC}.

	The process of finding a substitution rule for a tiling \tile\ is a two stage process: identifying candidates for the expansion predecessor and then examining which of these can be related to \tile\ by a replacement rule. An expansion predecessor tiling $\ol{\tile}$ must be locally isomorphic to $\lambda \tile$, where $\lambda$ is the expansion factor of the substitution rule. Moreover, as we consider vertex hierarchic substitution rules, any candidate for expansion predecessor must also satisfy the inclusion $vert(\ol{\tile}) \subset vert(\tile)$. 
	
	Having found a tiling $\ol{\tile}$ that satisfies these two conditions, the question is whether there exists a replacement rule relating $\ol{\tile}$ to $\tile$, or equivalently if each tile in $\ol{\tile}$ is covered by the same set of tiles in \tile. The latter step possibly requires assigining labels to distinguish between tiles with the same geometric shape (interval of a given length). Finally, it must be verified that the substitution rule on the labelled protoset admits a substitution tiling (infinite predecessors can be found). It is important to note that this last fact does not directly follow from the existence of a substitution rule relating $\tile$ to $\ol{\tile}$: it needs to be shown that $\ol{\tile}$ has a predecessor with respect to the same substitution rule, and more generally the existence of an infinite sequence of predecessors. 
	
	It may well happen that one identifies a substitution formulated on  too large a protoset, due to unnecessarily many labels for certain geometric tile shapes. This observation leads us to discuss a natural notion of \emph{equivalence} between substitution rules. We propose here the notion of \emph{patch equivalence}, under which two substitution rules are taken to be equivalent if for any number of iterations of the substitution rule on the tiles from the protoset, the sets of obtained patches (while ignoring the labels) are the same. In order to formalize this notion, let $\pi$ denote the function that extracts the geometric shape of a labelled prototile (interval with given length) or labelled patch (sequence of adjacent intervals). 
	\begin{definition}
		[Patch equivalent substitution rules] \label{defn_patch_equiv} Let $\tau$ and $\tau'$ be two labelled protosets and $\sigma$ and $\sigma'$ be substitution rules on respectively $\tau$ and $\tau'$. Then $\sigma$ and $\sigma'$ are \emph{patch equivalent} if for every prototile $T \in \tau$ and all $n \in \nat$ there exists a $T' \in \tau'$ such that $\pi(\sigma^n T )=\pi((\sigma')^n T')$, and similarly if for every prototile $T' \in \tau'$ and all $n \in \nat$ there exists a $T \in \tau$ such that $\pi(\sigma^n T )=\pi((\sigma')^n T')$. 
	\end{definition}
	We call a substitution rule \emph{patch minimal} if there exists no patch equivalent substitution rule with fewer prototiles. It is readily verified that within an entire set of patch equivalent substitution rules for a substitution tiling, there exists a unique patch minimal one. 
	
	\subsection{One-dimensional substitution tilings with an interval projection structure} \label{2to1_IPT} 
	
\begin{figure}
		[htp] $$ 
		\begin{xy}
			\xyMarkedImport{} 
			\xyMarkedPos{1}*\txt\labeltextstyle{$V$}
			\xyMarkedPos{2}*\txt\labeltextstyle{$W$}
			\xyMarkedPos{3}*\txt\labeltextstyle{$Y$}
			\xyMarkedPos{4}*\txt\labeltextstyle{$\Omega$}
		\end{xy}
		$$ \caption[Interval Projection Staircase]{Illustration of a tiling with interval projection structure and its associated {interval projection staircase}. The tiling is constructed by the projection to the subspace $V$ of a {\em slice} of the lattice $\integer^2$. Here $V$ here is the line with gradient $\frac{2}{1+\sqrt{5}}$ (the inverse of the golden ratio) and the slice is defined by the intersection of the lattice with a strip parallel to $V$. The lattice points within this strip are projected orthogonal, along the subspace $W=V^\perp$. By connecting lattice points within the strip whose projections yield neighbouring vertices in the projection tiling by line segments, we obtain the interval projection staircase.  The window is the subset of the space $Y$ that gives the slice when translated along $V$.} \label{2to1_interval_staircase} 
	\end{figure}
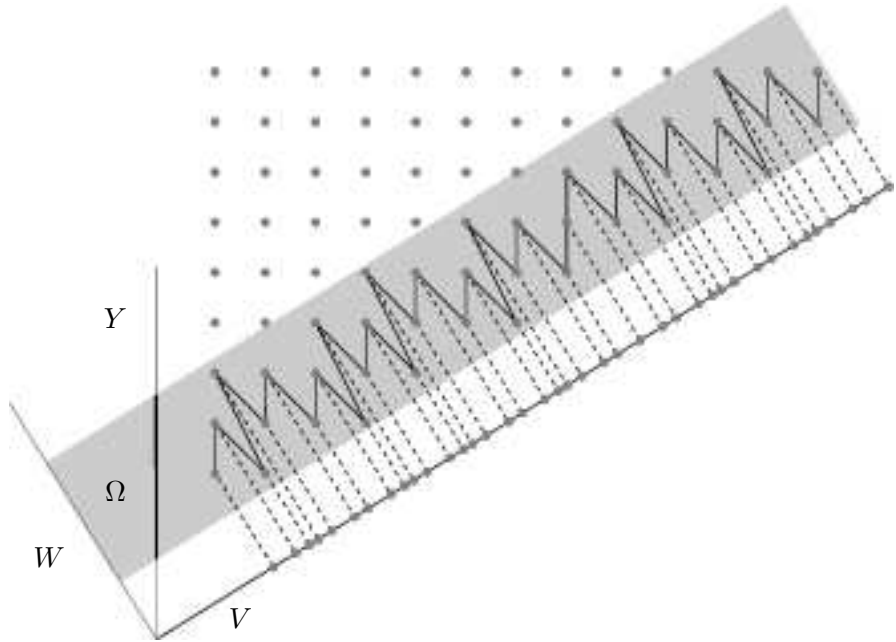

	Consider a one-dimensional tiling with an interval projection structure. The latter means that that the tiled line can be embedded in $\real^2$ in such a way that the vertices of the tiling are projections from points of the lattice $\integer^2$ that lie within a bounded distance of the embedded line. These lattice points may be linked together by line segments in $\real^2$ to form a staircase, as illustrated in Figure \ref{2to1_interval_staircase}. 
	
	\WARMprocessEPS{2to1_interval_staircase}{eps}{bb} 
	
\begin{figure}
		[htp] $$ 
		\begin{xy}
			\xyMarkedImport{} \xyMarkedPos{1}*!R\txt\labeltextstyle{Strip for $a$:} \xyMarkedPos{2}*!R\txt\labeltextstyle{Strip for $a+b$:} \xyMarkedPos{3}*!R\txt\labeltextstyle{Strip for $b$:} \xyMarkedPos{4}*!R\txt\labeltextstyle{Whole strip} 
		\end{xy}
		$$ \caption[Interval Projection Staircase window partition]{
The window partition (and associated partition of the strip) for the three possible staircase steps (to the right) of the interval projection staircase in Figure~\ref{2to1_interval_staircase}. These three types of steps give rise to three geometric tile types in the projection tiling. Note that the boundaries of the strip partition are obtained by translating the edges of the strip by two of the three staircase steps. The corrsponding window partition is obtained by projection to $Y$.} \label{2to1_window_eg} 
	\end{figure}
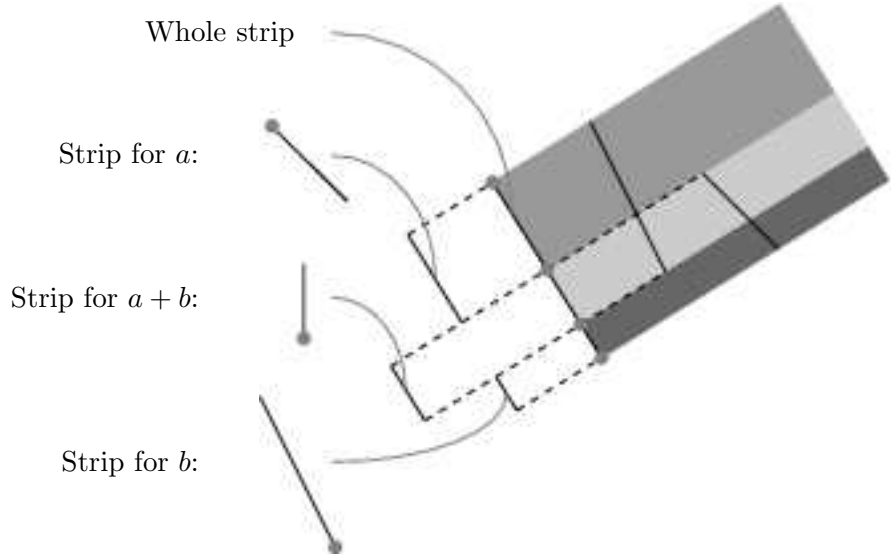

	The vertices of this staircase can be projected to the window space $Y$ that is transversal to $V$. The neighbourhood of a vertex within the staircase (and the neighbourhood of the corresponding vertex in the tiling) can be deduced entirely from the position of its projection along $V$ to $Y$. For instance, one can partition the window in regions indicating which staircase step lies to the right of a vertex whose projection to $Y$ lies in any of these regions. In Figure \ref{2to1_window_eg}, one finds such a window partition for the three possible right steps of the staircase in Figure~\ref{2to1_interval_staircase}. \WARMprocessEPS{2to1_window_eg}{eps}{bb} 
	
	The example in Figure~\ref{2to1_interval_staircase} displays three different staircase steps, projecting to three different interval tiles in $V$. It turns out that this number of staircase steps (and tile shapes) is in fact typical for one-dimensional tilings with an interval projection structure. This is the content of the next Lemma, which also features in Masakova {\em et al.} \cite{Masakova:LLPOQ,Masakova:SRFAS} and is related to the three distance theorems in number theory. For a historical account of the latter, see \cite[pp53--54]{Allouche:AS}. We include a short and self-contained proof for completeness. 

	\begin{lemma}
		\label{three_steps} A one-dimensional projection tiling with an interval projection structure has tiles of at most three lengths. In the case of three lengths, one is equal to the sum of the two others.
	\end{lemma}
	
	\begin{proof}
		Consider the points $\Pi_Y(t)$ and $\Pi_Y(t+l)$ at the ends of the interval window $\Omega_{l,t}\subset W$ with length $\Pi_Y(l)$, where $l\in\real^2$ and base $t\in\real^2$, and the translated lattice $\integer^2 + t$ (which has a point in common with the base of $\Omega_{l,t}$). Now consider the set of points $(\Omega_{l,t} + V) \cap (\integer^2+t)$, and its projection (along $W$) to $V$. 
 
		Since $\Omega_{l,t}$ is bounded, this projection yields a discrete point set in $V$. We now consider the two points in $\integer^2+t$ that project to either side of $\Pi_{V}(t)$ in $V$. We denote these points $t + a$ and $t - b$, where $a,b\in\integer^2$ correspond to steps of the staircase. 
  
		For simplicity, we continue the proof assuming that the window is the half-open interval $\Omega_{l,t}=[t,t+l)\subset W$. The argument for the remaining case in which $\Omega_{l,t}=(t,t+l]\subset W$ is analogous. 
		
		The subwindows corresponding to these two steps are the intervals $\Omega_{l,t}^a\subset \Omega_{l,t}$ and $\Omega_{l,t}^b\subset \Omega_{l,t}$ given by 
		\begin{eqnarray*}
			\Omega_{l,t}^a&=&\{p\in\Omega_{l,t}~|~\Pi_{Y}(a)+p\in\Omega_{l,t}\}=[\Pi_Y(t), \Pi_{Y}(t+l -a)) ,\\
			\Omega_{l,t}^b&=&\{p\in\Omega_{l,t}~|~\Pi_{Y}(b)+p\in\Omega_{l,t}\}=[ \Pi_{Y}(t - b), \Pi_Y(t+l)). 
		\end{eqnarray*} 
		If $\Pi_{Y}(t+l - a) = \Pi_Y(t -b)$ then all steps in the staircase are accounted for, and there are only two types of steps in the staircase. \WARMprocessEPS{2to1_three_steps_window}{eps}{bb} 
		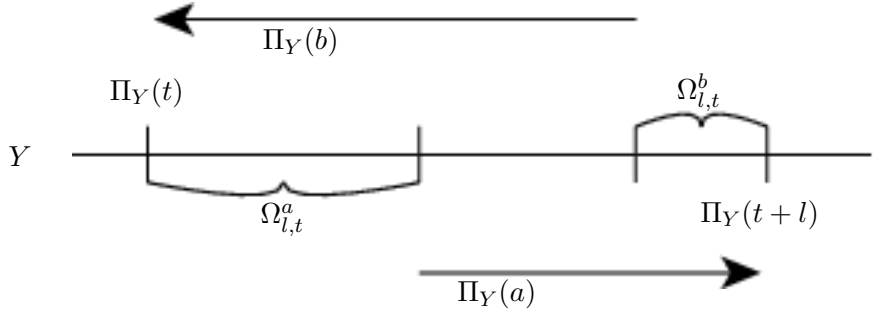
\begin{figure}
			[htp] $$ 
			\begin{xy}
				\xyMarkedImport{} \xyMarkedPos{1}*!L\txt\labeltextstyle{$Y$} \xyMarkedPos{2}*\txt\labeltextstyle{$\Pi_Y(t)$} \xyMarkedPos{3}*\txt\labeltextstyle{$\Pi_Y(t+l)$} \xyMarkedPos{4}*\txt\labeltextstyle{$\Pi_Y(a)$} \xyMarkedPos{5}*\txt\labeltextstyle{$\Pi_Y(b)$} \xyMarkedPos{8}*\txt\labeltextstyle{$\Omega^a_{l,t}$} \xyMarkedPos{9}*\txt\labeltextstyle{$\Omega^b_{l,t}$} 
			\end{xy}
			$$ \caption[Three steps window]{The partition of the window for a staircase with three steps. The arrows indicate how nearest neighbour staircase steps project to the window. The central subwindow is the subwindow for staircase step $a + b$.} \label{2to1_three_steps_window} 
		\end{figure}
		
		If $\Pi_{Y}(t+l - a) \neq \Pi_Y(t -b)$, the situation is as shown in Figure~\ref{2to1_three_steps_window}. Apart from the subwindows for $a$ and $b$ discussed above, there is a remaining third subwindow in between. Consider a point $q\in\real^2$ whose projection to $Y$ lies in this central subwindow, i.e. $\Pi_Y(q) \in [\Pi_Y(t+l - a) , \Pi_Y(t - b))$. Then the projection to $Y$ of the point $q + a + b$ lies in the interval $[\Pi_Y(t+l + b) , \Pi_Y(t + a))$ which is a subinterval of $\Omega_{l,t}$. We finally claim that if $q\in\integer^2$ projects to the central subwindow, the next (right) step in the staircase is $a+b$. Namely, if there are points whose projection to $V$ lies between $\Pi_V(q)$ and $\Pi_V(q + a + b)$, then by periodicity of the lattice $t + a$ and $t - b$ would not project to the closest points to $t$, contradicting the definition of $a$ and $b$. 
		
		Hence, there are three types of staircase steps: $a$, $b$ and $a + b$. The corresponding subwindows are $\Omega_{l,t}^a=[\Pi_Y(t), \Pi_Y(t+l - a))$, $\Omega_{l,t}^{a+b}=[\Pi_Y(t+l - a),\Pi_Y(t - b))$ and $\Omega_{l,t}^b= [\Pi_Y(t - b), \Pi_Y(t+l))$. 
		
		In the case of the other choice of the half-open window interval $\Omega_{l,t}=(\Pi_Y(t),\Pi_Y(t+l)]$, one finds $\Omega_{l,t}^a=(\Pi_Y(t), \Pi_Y(t+l - a)]$, $\Omega_{l,t}^{a+b}=(\Pi_Y(t+l - a),\Pi_Y(t - b)]$ and $\Omega_{l,t}^b= (\Pi_Y(t - b), \Pi_Y(t+l)]$. 
	\end{proof}
	
	The staircase can be seen as the orbit of a function $\Phi_{l,t}$, which maps a vertex of the staircase to its right nearest neighbour. 

	\begin{definition}[The staircase function $\Phi_{l,t}$] \label{defn:Phi}
		
		The function $\Phi_{l,t}: (V + \Omega_{l,t}) \to (V + \Omega_{l,t})$ is defined as $$ 
		\begin{array}
			{lcc} \Phi_{l,t}(z)&=& \left\{ 
			\begin{array}
				{ll} z + a & \text{if }\Pi_Y(z) \in \Omega_{l,t}^a,\\
				z + a + b & \text{if }\Pi_Y(z) \in\Omega_{l,t}^{a+b},\\
				z + b & \text{if }\Pi_Y(z) \in \Omega_{l,t}^b
			\end{array}
			\right. 
		\end{array}
		$$ where $\Omega_{l,t}^a$, $\Omega_{l,t}^{a+b}$, and $\Omega_{l,t}^b$ are the subwindows of $\Omega_{l,t}$ for the steps $a$, $a+b$ and $b$, as introduced in the proof of Lemma~\ref{three_steps} above.
	\end{definition}
	
	As the function $\Phi_{l,t}$ is invertible, by iterating $\Phi_{l,t}^{-1}$ one can also generate successions of neighbours to the left. $\Phi_{l,t}$ is defined for every point $r$ inside the strip $V + \Omega_{l,t}$, so that the set $\{\Phi_{l,t}^n(r)~|~n\in\integer\}$ contains the vertices of an interval projection staircase in $r + \integer^2$. 

	This is just a translation of the interval projection staircase in $\integer^2$ with window $\Omega_{l,t-\Pi_Y(r)}$. 

	The subwindows for $a$, $b$ and $a+b$ can be further subdivided, providing a partition relating projections of vertices of the staircase to patches (collections of adjacent tiles) instead of single tiles (or staircase steps). 
 
	Let us illustrate this in more detail. 

	Consider the interval $[\Pi_Y(t), \Pi_Y(t+l - a))$ in $Y$. The vertices that project to this subwindow are followed by the step $a$. Call the corresponding set of vertices of the staircase $S_a$. Because of the definition of $S_a$ we have $\Phi_{l,t}(S_a)=S_a + a$, representing the set of vertices at the end of an $a$ step. The projection $\Pi_Y(\Phi_{l,t}(S_a))$ densely fills the window interval $[\Pi_Y(t+a), \Pi_Y(t+l))$, the subwindow for $a$ translated by $\Pi_Y(a)$. 
	
	For each of the vertices in $S_a+a$, the next right step (after $a$) depends on the subwindow to which it projects. If one or both of the points $\Pi_Y(t+l - a)$ and $\Pi_Y(t - b)$, lie in $[\Pi_Y(t+a), \Pi_Y(t+l))$ then the latter interval can be partitioned into several subintervals, on the basis of their projection to one of the three subwindows for $a$, $b$ and $a+b$. 
The image of these subintervals under $\Phi_{l,t}^{-1}$ yields a partition of the subwindow for $a$, signalling the next two steps that follow a vertex (the first one being $a$). 
  
	This process can be continued by induction to construct window partitions, corresponding to sequences of steps of the staircase (and corresponding tiles in the projection tiling) after a given vertex. The points in the window that act as boundary points for these partition intervals are the images of the two boundary points under iterates of the map $\Phi_{l,t}^{-1}$: $\Pi_Y(\Phi_{l,t}^{-n}(\Pi_Y(t+l - a))$ and $\Pi_Y(\Phi_{l,t}^{-n}(\Pi_Y(t - b))$. We may therefore associate tiling (and staircase) patches with window intervals. 
	\begin{definition}
		[Window Interval for the Patch $P$] The \emph{window interval} $I_P$ for the patch $P$ of a one-dimensional tiling with interval projection structure is the unique interval of the window such that $\Pi_Y(z) \in I_P$ for some vertex $z$ of the projection staircase, if and only if the vertex $\Pi_V(z)$ of the tiling is followed (to the right) by the patch $P$. 
	\end{definition}
	We now use the window intervals to show that the set of one-dimensional tilings with an interval projection structure with fixed $V$, $W$ and $l$ forms a local isomorphism class. 
	\begin{lemma}
		\label{2to1_LI_class} The set of one-dimensional interval projection tilings characterised by the linear subspaces $V$ and $W$ of $\real^2$ satisfying $\real^2=V\oplus W$, and fixed $l \in \real^2$, forms a local isomorphism class of one-dimensional tilings. 
	\end{lemma}
	\begin{proof}
		Let $\tile_{l,t}$ denote the interval projection tiling with window vector $l$ and base $t$. We first show that the tilings $\tile_{l,t}$ and $\tile_{l,t'}$ are locally isomorphic for any $t, t' \in \real^2$. Consider a patch $P$ of $\tile_{l,t}$, with window interval $I_P \subset \Omega_{l,t}$. The equivalent window interval $\Pi_Y(I_P - t + t')$ can be found in $\Omega_{l,t'}$. Any vertex of the staircase for $\tile_{l,t'}$ which projects to this interval will be followed by the patch $P$, and as $\Pi_Y(\integer^2)$ is dense, such vertices exist. Thus every patch of $\tile_{l,t}$ occurs in $\tile_{l,t'}$, by an analogous argument any patch of $\tile_{l,t'}$ occurs in $\tile_{l,t}$ so the two tilings are locally isomorphic. 
		
		We now show that any tiling of the local isomorphism class has an interval projection structure. Consider a tiling \tile that is locally isomorphic to $\tile_{l,t}$. The tiling must have the same tile shapes so we can find a staircase $\mathcal{S}$, with the same steps as the staircase for $\tile_{l,t}$, such that $\Pi_V(\mathcal{S}) = \tile$. Now consider the closure of the projection of $\mathcal{S}$ to $Y$: $\ol{\Pi_Y(\mathcal{S})} \subset W$. If the extreme points of this set were more than $\Pi_Y(l)$ apart, then there would be two vertices of $\mathcal{S}$ laying more than $\Pi_Y(l)$ apart in the projection to $Y$. This is impossible, as these points would lie at the ends of some finite patch, which would appear in $\tile_{l,t}$. We need only show, therefore, that the set $\ol{\Pi_Y(\mathcal{S})}$ is an interval. To do this consider the projection of longer and longer finite lengths of the staircase to $Y$. These have smaller and smaller gaps between the points as each finite length of staircase is also part of the staircase for $\tile_{l,t}$, whose window is an interval. We may therefore get arbitrarily close to any point in the interval between the extreme points of $\ol{\Pi_Y(\mathcal{S})}$ so that the closure yields the interval. 
		
		This concludes the proof of the fact that all tilings in $\{\tile_{l,t}|t \in \real^2\}$ are locally isomorphic and any tiling locally isomorphic to $\tile_{l,t}$ is contained in the set $\{\tile_{l,t}|t \in \real^2\}$.

	\end{proof}
	Note that if $\Pi_Y(t) = \Pi_Y(t'+z)$ for some $z \in \integer^2$ then $\tile_{l,t}$ and $\tile_{l,t'}$ differ only by a translation. Thus, with respect to the natural equivalence relation (see footnote~\ref{fn:equivalence}) each value of $t\in Y/\Pi_Y(\integer^2)$ is associated with a unique tiling. 

\section{The partition algorithm} \label{partition}
	
In this section we discuss substitution rules in the context of the geometry of the interval projection structure. Section \ref{partition_defn} provides a discussion of the setting, after which in Section~\ref{partition_alg}, we describe the existence of substitution rules, in this context. In Section~\ref{sec:proof_subst} we examine precisely when substitutions rules exist, and provide a chacterisation of all substitution rules for any given tiling of the line with an interval projection structure.
	
The proofs of our results are all constructive, and in Section~\ref{partition_eg} we illustrate the construction of a substitution rule in the case of a specific tiling. The algorithm has been implemented in Mathematica and is available from the authors \cite{Harriss:NSTOT_code}.
	
	\subsection{Window partitions and replacement rules} \label{partition_defn}
	
	In this section we show that every one dimensional tiling with a interval projection structure and substitution rule can be associated to a unimodular $2\times2$ primitive integer matrix. The action of this matrix on the projection structure is key to providing a geometric interpretation of the substitution rule. 
	
	We first associate a matrix to any substitution tiling with an interval projection structure. This procedure is folklore, see for instance \cite{Kenyon:SST}. 

	\begin{lemma}
		\label{2to1_matrix_lemma} The spaces $V$ and $W$ of a nonperiodic\footnote{Recall that nonperiodic refers to the geometry of the tiling before adding any labelling.} one-dimensional projection tiling with interval projection structure with a substitution rule are the eigenspaces of a non-singular primitive, $2\times2$ integer matrix. 
	\end{lemma}
	\begin{proof}
		
	Consider a tiling $\mathcal{T}$ with interval projection structure from the lattice $\mathcal{L}$ and substitution rule $\sigma$. Let $\mathcal{S}$ be the staircase for $\mathcal{T}$. Now consider the expansion predecesssor $\ol{\mathcal{T}}$ for $\mathcal{T}$ under $\sigma$. We may lift $\ol{\mathcal{T}}$ to a staircase $\ol{\mathcal{S}}$ whose vertices form a subset of the vertices of $\mathcal{S}$. This staircase lies on a lattice $\ol{\mathcal{L}}$ which is a sublattice of $\mathcal{L}$. The steps in this lattice corresponding to the tiles $a$ and $b$ in $\mathcal{T}$ generate this lattice, just as the equivalent steps in $\mathcal{S}$ generate $\mathcal{L}$. We may therefore define the linear map taking $\mathcal{L}$ to $\ol{\mathcal{L}}$ as the map taking the steps corresponding to tiles $a$ and $b$ in $\mathcal{S}$ to the equivalent steps in $\ol{\mathcal{S}}$.
		
		The action of this matrix must commute with the projection, thus the spaces $V$ and $W$ are the eigenspaces of the matrix. As the tiling is non-periodic, the space $V$ must be irrational, and the matrix must therefore be primitive and non-singular. 
	\end{proof}
	
	A key result of this paper concerns the fact that the matrix mentioned in the previous lemma must be unimodular, and thus represents an automorphism of the lattice in the embedding space. In previous studies \cite{Pleasants:DQ, Masakova:SRFAS, Luck:TNOTA} unimodularity of the matrix was a hypothesis. 

	A substitution rule consists of two components: inflation and replacement. First the tiling is homogeneously inflated by a constant factor. In the embedding space we may represent this inflation as multiplication by a matrix $M$, as identified in Lemma~\ref{2to1_matrix_lemma}.  This maps the staircase $\mathcal{S}$ for the tiling $\mathcal{T}$ to $M \mathcal{S}$ which serves as a staircase for the inflated tiling. 

In the embedding space the replacement rule corresponds to replacing every staircase step in $M \mathcal{S}$ by a section of staircase steps with vertices on the lattice $\lattice$.
	
	We now discuss how the substitution rule acts on the window.
	
	The existence of a substitution rule for an interval projection tiling implies the existence of a patch minimal substitution rule for the same tiling. So we assume without loss of generality that $\sigma$ is patch minimal: for some iteration $\sigma^n$ of $\sigma$ the prototiles can be distinguished by the patches they substitute to. Each of these patches is associated to a particular subwindow of $\Omega$ (where $\Omega$ is the window for $\mathcal{T}$). 
	
	Let $\ol{\mathcal{T}}$ be the expansion predecessor of $\mathcal{T}$ with respect to substitution rule $\sigma^n$. Every instance of a prototile $\ol{T}$ in $\ol{\mathcal{T}}$ is replaced by the same patch $P$ of tiles in $\mathcal{T}$. However $P$ can also arrise at other places in the tiling where it is not covered by $\ol{T}$. For example $P$ could occur as a subset of a patch replacing a different prototile. This occurs in the Fibonacci tiling with the substitution rule $(a\to ab,b\to a)$: the tile $\ol{b}$ in the expansion predecessor is replaced by the tile $a$, but $a$ is also contained in the patch $ab$ that replaces $\ol{a}$. Thus the subwindow $\Omega_P$ for $P$ within $\Omega$ is not necessarily the subwindow $\ol{\Omega_T}$ within $\ol{\Omega}$ associated to the prototile $\ol{T}$, but it must contain it. 
	
	By Lemma \ref{three_steps} the window $\ol{\Omega}$ may be divided up into two or three subwindows corresponding to the lengths of the tiles. Let $\ol{\Omega_{|T|}}$ be the subwindow for all tiles of the same length as $\ol{T}$. The intersection $\ol{\Omega_{|T|}} \cap \Omega_P$ corresponds to every tile of the same length as $\ol{T}$ in $\ol{\mathcal{T}}$ that is covered by $P$. Consequently, we have $\ol{\Omega_T} = \ol{\Omega_{|T|}} \cap \Omega_P$.
	
	The replacement rule takes the instances of $T$ and replaces each of them with the same patch of tiles. In the staircase the corresponding steps in $\ol{\mathcal{S}}$ are each replaced by a collection of steps $\{s_1,\ldots,s_n\}$ in $\mathcal{S}$. The vertices of the $s_i$'s project to $Y$ (after closure) as the subintervals $\Pi_Y(s_i) + \ol{\Omega_T}$ of the window $\Omega$. Furthermore, as $\sigma^n$ is a substitution rule, each $\Pi_Y(s_i) + \ol{\Omega_T}$ is contained within the interval associated to a unique prototile of $\mathcal{T}$.

        \begin{theorem}
                \label{2to1_matrix} The spaces $V$ and $W$ of a nonperiodic one-
dimensional projection tiling with interval projection structure with  
a substitution rule are the eigenspaces of a primitive matrix in $Gl
(2,\integer)$.
        \end{theorem}
        \begin{proof} 
From Lemma \ref{2to1_matrix_lemma} we have a matrix $M\in gl(2,\integer)$, it remains to be shown that $M$ is unimodular.
                
                The vertices of the staircase for the tiling $\tile$ are the points  
in $\integer^2$ which project to the window $\Omega$ in $Y$. The  
vertices of the staircase for the expansion predecessor $\overline
{\tile}$ are the vertices of $M \integer^2$ which project to the  
subwindow $\ol{\Omega}$. The lattice $M \integer^2$ is a subgroup of $
\integer^2$ and both are Abelian groups. We may therefore consider $
\integer^2$ as a finite union of distinct translations of $M \integer^2$, 
$$\integer^2=M  
\integer^2 \cup (g_1 + M \integer^2) \cup\ldots\cup (g_{N-1} + M  
\integer^2).$$ We denote the corresponding set of translations as
$G=\{0,g_1,\ldots,g_{N-1}\}$. Note that $N=|\det(M)|$.
                
The replacement rule acting on $\ol{\tile}$ yields $\tile$.  By the  
discussion above we may view the replacement procedure in the 
projection to the window $\Omega \subset Y$.  The vertices of a labelled
predecessor staircase step project to the points $\Pi_Y(M \integer^2) \cap [a,b] \subset \ol{\Omega}$.
The replacement rule induces translations of this interval by elements of
$\Pi_Y(\integer^2)$.  Consider one such translation, by $\Pi_Y(z)$, for some $z$ in $\integer^2$.  We 
may associate $g_i \in G$ to $z$ where $z \in M \integer^2 + g_i$.  
This translation adds to the vertices of the staircase 
the points in $M \integer^2 + g_i$ that project to $[a,b] + Pi_Y(z)$ in $Y$.
In order to get every point in $\Pi_Y(\integer^2)$ that projects to $\Omega$,  
therefore, the action of the replacement rule on the subwindow must  
cover the window at least $N$ times.  Furthermore the set of  
intervals covering any region must include an interval translated by  
an element of $\Pi_Y(M \integer^2 + g_i)$ for each $g_i \in G$.

Let $\mathcal{P}$ be the set of intervals induced by the action of replacement  
rule on the labelling intervals of the predecessor.  For $g \in G$, let the set of intervals in $\mathcal{P}$  
translated by $z \in \Pi_Y(M \integer^2 + g)$ be $\mathcal{P}_g$.  
The labelling partition of the window $\Omega$ is defined by a finite number of  
points $\{v_1,\ldots,v_p\}$.  The subwindow $\ol{\Omega}$ is  
partitioned by the points $\{\ol{v}_1,\ldots,\ol{v}_p\}$ in the same  
manner.  Each translate of a labelling interval in the subwindow partition  
induced by the replacement rule corresponds to a unique labelled tile  
and is thus the subset of an interval of the partition on $\Omega$.  
Consider $v_1$.  As the end point of an interval of the partition on $
\Omega$, it must be the end point of any interval in $\mathcal{P}$  
that contains it.  Consider the intervals for which $v_1$ is a left  
end point.  These must contain an element of $\mathcal{P}_g$ for each  
$g \in G$.  The left end points of these intervals are the  
translations of points in $\Pi_Y(\{\ol{v}_1,\ldots,\ol{v}_p\})$.  
Thus for each $g \in G$ there exists $\ol{v}_g \in \{\ol{v}_1,\ldots,
\ol{v}_p\}$ and $z_g \in M \integer^2$ such that $v_1 = \ol{v}_g +  
\Pi_Y(z_g + g)$.  Furthermore, $\ol{v}_g$ is equal to $\ol{v}_h$ if  
and only if $g=h$, as elements of $G$ are not in $M \integer^2$.

	The subwindow associated to the $q$-th iterate of the substitution  
rule is partitioned in the same manner as $\Omega$ by $p$ points. The  
absolute value of the determinant of the matrix for this substitution  
rule is $N^q$. Thus the lattice is partitioned into $N^q$ translates  
of $M^q \integer^2$. Thus if $N > 1$, there exists $q$ such that $N^q  
> p$.  This is leads to a contradiction, as $\ol{v}_g$ is equal to $\ol{v}_h$ if  
and only if $g=h$, so every element of $g$ must have a distinct $\ol
{v}_g$.  Hence $N = 1$ and consequently $M$ is unimodular.
\end{proof}
	
	The following lemma describes the set of potential expansion predecessors: tilings, whose vertices are a subset of the original tiling, and locally isomorphic up to scaling. 
	\begin{lemma}
		\label{2to1_expansion} Let $M$ be a primitive matrix in $Gl(2,\integer)$, with eigenspaces $V$ and $W$ and expanding eigenvalue $\lambda$. Let $\Omega_{l,t}$ be the window of the interval projection tiling $\tile_{l,t}$ with window and tiling spaces $V$ and $W$ and $t \in \real^2$. Let $\ol{\Omega}_{l,t,s} = M \Omega_{l,0} + \Pi_Y(s + t)$, be the window of the interval projection tiling $\ol{\tile}_{l,t,s}$, where $s \in \real^2$ such that $\ol{\Omega}_{l,t,s} \subset \Omega_{l,t}$. The interval projection tilings $\lambda \tile_{l,t}$ and $\ol{\tile}_{l,t,s}$ are locally isomorphic for all $t \in \real^2$. Furthermore $vert(\ol{\tile}_{l,t,s}) \subset vert(\tile_{l,t})$. 
	\end{lemma}
	\begin{proof}
		By Lemma \ref{2to1_LI_class} the tilings $\{\tile_{l,t}|t \in \real^2\}$ are locally isomorphic.
		
		We show that every tiling in $\{\ol{\tile}_{l,t,s}| t \in \real^2, s \in \real^2\}$ is contained in $\{\lambda \tile_{l,t}|t \in \real^2\}$. Let $\mathcal{S}_{l,t}$ be the staircase associated to $\tile_{l,t}$. As the action of $M$ commutes with the projection, we have $\Pi_V(M \mathcal{S}_{l,t}) = \lambda \tile_{l,t}$ and $\ol{\Pi_Y(M \mathcal{S}_{l,t})} = \pm \lambda^{-1} \Omega_{l,t}$. Hence $$ \ol{\Omega}_{l,t,s} = \Pi_Y M \Omega_{l,0} + \Pi_Y(s) + \Pi_Y(t) = \Pi_Y M \Omega_{l,\Pi_Y(s) + t} $$ Thus $\ol{\tile}_{l,t,s} = \lambda \tile_{l,s+t}$, so every tiling in $\{\ol{\tile}_{l,t,s}| t \in \real^2, s \in \real^2\}$ is a tiling in $\{\lambda \tile_{l,t}|t \in \real^2\}$.
		
		Finally, if $\ol{\Omega}_{l,t,s} \subset \Omega_{l,t}$ then every vertex of $\ol{\tile}_{l,t,s}$ is a vertex of $\tile_{l,t}$. 
	\end{proof}

	\subsection{The algorithm} \label{partition_alg}
	
	Lemma \ref{2to1_expansion} describes the set of potential expansion predecessor tilings. It remains to show for which of these there exists a replacement rule. Consider the window $\Omega$ with subwindow $\ol{\Omega}$, for a potential expansion predecessor. Both $\ol{\Omega}$ and $\Omega$ are partitioned into two or three intervals associated to the two or three tile types given by Lemma \ref{three_steps}. The tiles of $\ol{\mathcal{T}}$ are covered by patches of $\mathcal{T}$. These are finite patches of $\mathcal{T}$ so they have interval windows. The subwindow $\ol{\Omega}$ is therefore partitioned into intervals each of which is associated to a specific tile type in the expansion predecessor with a specific patch of tiles that cover it. We use this partition to label the tiles in $\ol{\mathcal{T}}$. The identification of an appropriate rescaling of $\Omega$ with $\ol{\Omega}$ induces a corresponding labelling partition of $\Omega$. In turn this provides a labeling of the tiles in $\mathcal{T}$. 
	
	In the procedure above we initially identify tile types only by length. We now repeat the procedure, initially identifying tile types by length and label. This produces a labelling of the tiles, that is either the same as the previous labelling or a refinement of it. In the former case we have identified a labelling with respect to which the tiling admits a substitution rule. In the latter case we may repeat the procedure. The tiling admits a substitution rule for this potential expansion predecessor if and only if after a finite number of iterations no further refinement of labelling occurs. 
	
	To describe this procedure in detail, we reintroduce the points $l$, $t$ and $s$ that describe the position of the window $\Omega_{l,t}$ and the subwindow $\ol{\Omega}_{l,t,s}$. We define a function $g_{l,t,s}: V + \Omega_{l,t} \to V + \Omega_{l,t}$ by: $$ g_{l,t,s}(z) = M^{-1}(\Phi_{l,t}^{-n}(z) - s - t) + t, $$ where $$ n = \min\{n \in \integer^+ | \Pi_Y(\Phi_{l,t}^{-n}(z)) \in \ol{\Omega}_{l,t,s}\}. $$ and $\Phi_{l,t}^{-1}$ is the staircase function defined in Definition \ref{defn:Phi}.
	
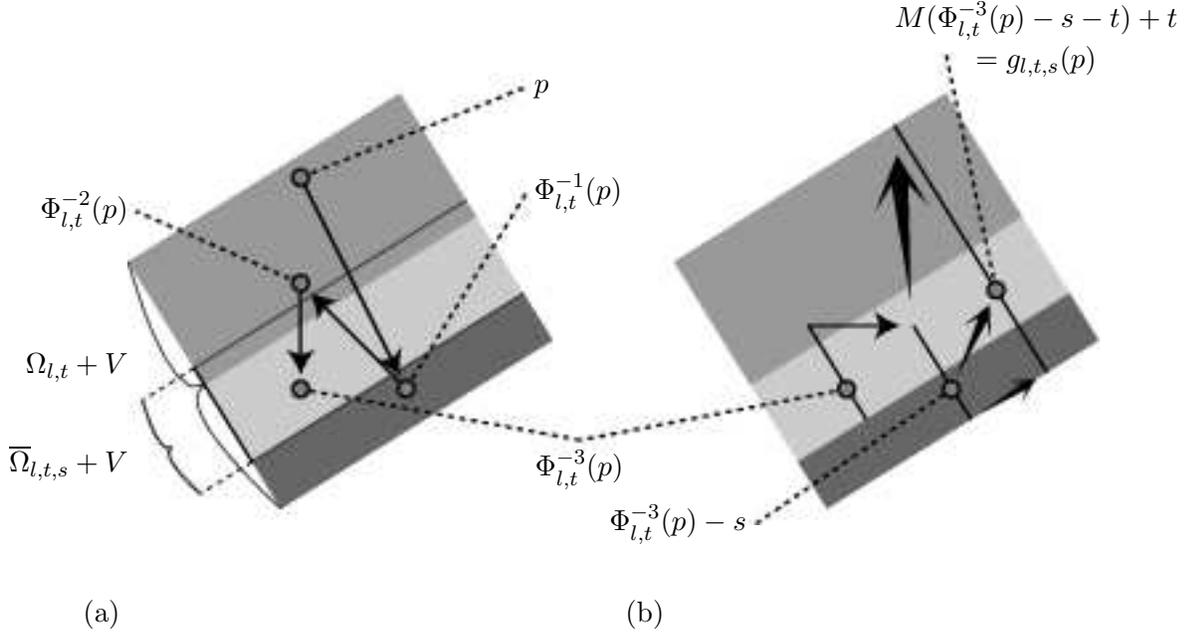
\begin{figure}
		[htp] $$ 
		\begin{xy}
			\xyMarkedImport{} 
			\xyMarkedPos{1}*+!L\txt\labeltextstyle{$p$}
			\xyMarkedPos{2}*+!L\txt\labeltextstyle{$\Phi_{l,t}^{-1}(p)$} 
			\xyMarkedPos{3}*+!R\txt\labeltextstyle{$\Phi_{l,t}^{-2}(p)$} 
			\xyMarkedPos{4}*+!U\txt\labeltextstyle{$\Phi_{l,t}^{-3}(p)$} 
			\xyMarkedPos{5}*+!R\txt\labeltextstyle{$\Phi_{l,t}^{-3}(p)-s$} 
			\xyMarkedPos{6}*+!L\txt\labeltextstyle{$M (\Phi_{l,t}^{-3}(p) - s - t) + t$ \\
			$= g_{l,t,s}(p)$} \xyMarkedPos{8}*+!R\txt\labeltextstyle{$\Omega_{l,t} + V$} 
			\xyMarkedPos{7}*+!R\txt\labeltextstyle{$\ol{\Omega}_{l,t,s} + V$} 
			\xyMarkedPos{10}*+!R\txt\labeltextstyle{(a)} 
			\xyMarkedPos{11}*+!R\txt\labeltextstyle{(b)} 
		\end{xy}
		$$ \caption[The function $g_{l,t,s}$]{Illustration of the action of the function $g_{l,t,s}$ to a point $x \in V + \Omega_{l,t}$. (a) We consider a point $p$ in the $V+\Omega_{l,t}$ and then apply $\Phi_{l,t}^{-1}$ to walk the point back until the image lies in $\ol{\Omega}_{l,t,s} + V$. Here $\Phi_{l,t}^{-1}(p)$ and $\Phi_{l,t}^{-2}(p)$ do not lie $\ol{\Omega}_{l,t,s} +V$ but $\Phi_{l,t}^{-3}(p)$ does. (b) Subsequent application of the affine transformation.} \label{2to1_applying_g} 
	\end{figure}

	To illustrate the action of this function, consider $x \in \integer^2 \cap V + \Omega_{l,t}$. Using the staircase function $\Phi_{l,t}^{-1}$, we follow points on the staircase until we reach a point $p$ that projects to the subwindow $\ol{\Omega}_{l,t,s}$. We apply an affine transformation, based on $M^{-1}$, that takes $V + \ol{\Omega}_{l,t,s}$ to $V + \Omega_{l,t}$. This affine transformation relates points in $\Omega_{l,t} + V$ and $\ol{\Omega}_{l,t,s} + V$, whose $\Pi_Y$-projections have the same relative position in their respective windows. The action of the function $g_{l,t,s}$ is illustrated in Figure \ref{2to1_applying_g}. 
	
	\WARMprocessEPS{2to1_applying_g}{eps}{bb} 
	
	Consider the points $t+l-a$ and $t-b$, whose $\Pi_Y$-projections give the boundary points of the partition of the window $\Omega_{l,t}$ associated to tile type, as in Lemma \ref{three_steps}. Let $G_{l,t,s}$ denote the subset of $\Omega_{l,t}$ to which all $g_{l,t,s}$ iterates project, $$ G_{l,t,s} = \Pi_Y(\{g_{l,t,s}^n(t+l - a)| n \in \integer^+\}) \cup \Pi_Y(\{g_{l,t,s}^n(t - b)| n \in \nat\}). $$
	
	As, $G_{l,t,s} = G_{l,0,s} + \Pi_Y(t)$, the relative position of $G_{l,t,s}$ as a subset of $\Omega_{l,t}$ is independant of $t$. 
	
	Finiteness of $G_{l,t,s}$ is a necessary and sufficient condition for the existence of a substitution rule. 
	\begin{lemma}
		\label{2to1_finiteifsubst} The set $G_{l,t,s}$ is finite if and only if there exists a substitution rule for the interval projection tiling $\mathcal{T}_{l,t}$, with expansion predecessor $\ol{\mathcal{T}}_{l,t,s}$ associated to the subwindow $\ol{\Omega}_{l,t,s} \subset \Omega_{l,t}$. 
	\end{lemma}
	\begin{proof}
		We first prove that there is a substitution rule if the set $G_{l,t,s}$ is finite.
		
		If $G_{l,t,s}$ is finite, it provides the boundary points of a partition of $\Omega_{l,t}$ into a finite number of intervals. Consider such an interval $(g_1,g_2)$. In the subwindow $\ol{\Omega}_{l,t,s}$ the equivalent interval is $(\ol{g_1},\ol{g_2}) = M ((g_1,g_2) - t) + \Pi_Y(s + t)$. 
		
		We show that $(\ol{g_1},\ol{g_2}) \cap G_{l,t,s}$ is empty. Note that $$ g_{l,t,s}((\ol{g_1},\ol{g_2})) = (g1,g2) $$ as $(\ol{g_1},\ol{g_2}) \subset \ol{\Omega}_{l,t,s}$. Therefore $$ g_{l,t,s}((\ol{g_1},\ol{g_2}) \cap G_{l,t,s}) \subset ((g1,g2) \cap G_{l,t,s}). $$ Since $((g1,g2) \cap G_{l,t,s})$ is empty by construction, it follows that $(\ol{g_1},\ol{g_2}) \cap G_{l,t,s}$ is empty.
		
		As $\Pi_Y(t+l - a)$ and $\Pi_Y(t - b)$ are in $G_{l,t,s}$, and $(\ol{g_1},\ol{g_2}) \cap G_{l,t,s}$ is empty, $(\ol{g_1},\ol{g_2})$ must lie entirely in the subwindow associated to one of the tile types defined in Lemma \ref{three_steps}. Denote the associated staircase step $c \in \integer^2$. Consider the translation of $(\ol{g_1},\ol{g_2})$ by $\Pi_Y(c)$ $(\Pi_Y(\ol{g_1}+c),\Pi_Y(\ol{g_2}+c))$. We show that $(\Pi_Y(\ol{g_1}+c),\Pi_Y(\ol{g_2}+c)) \cap G_{l,t,s}$ is empty. Namely, suppose that this set contains a point $p$. Under $g_{l,t,s}$, $p$ is taken to $p-c$ by $\Phi_{l,t}^{-1}$, which is in the subwindow $\ol{\Omega}_{l,t,s}$, and then mapped to $M^{-1}(p - c - t) + s + t$ in $\Omega_{l,t}$. This yields a point in $((g1,g2) \cap G_{l,t,s})$, which is empty by construction.
		
		By induction of the above argument, there exists an $N$ such that for all positive $n < N$, $\Pi_Y \circ \Phi_{l,t}^n((\ol{g_1},\ol{g_2})) \cap G_{l,t,s}$ is empty, and $\Pi_Y \circ \Phi_{l,t}^N((\ol{g_1},\ol{g_2})) \subset \ol{\Omega}_{l,t,s}$. This implies that the partition of $\Omega_{l,t}$ with boundary points in $G_{l,t,s}$ is the set of labelling intervals for a substitution rule. 
		
		We now prove the converse: if $G_{l,t,s}$ is infinite there is no substitution rule.
		
		Suppose there is a substitution rule $\sigma$, with an associated labelling partition of the window $\Omega_{l,t}$. Let the vertices of this partition be the set of points $P$. The set $P$ must contain the points $\Pi_Y(t+l - a)$ and $\Pi_Y(t - b)$. The set $G_{l,t,s}$ is the union of the $g_{l,t,s}$ orbits of $\Pi_Y(t+l - a)$ and $\Pi_Y(t - b)$, one of which must be infinite. As $P$ is finite there must exist some point $p \in P$ such that $g_{l,t,s}(p) \notin P$. Consider the labelling interval that contains $g_{l,t,s}(p)$. The corresponding interval $I$ in the subwindow $\ol{\Omega}_{l,t,s}$ is translated to a series of intervals of the window $\Omega_{l,t}$, by the replacement rule for $\sigma$ (as discussed in more detail above). One of these translated intervals, contains $p$ in its interior. As $p \in P$, it is on the boundary between two labelling intervals. Thus we can find two tiles whose label is associated to $I$ that are covered by two different labelled patches after applying the replacement rule. This contradicts the assumption that $\sigma$ is a substitution rule. 
	\end{proof}

	\subsection{Illustration} \label{partition_eg}
	
	In this section we illustrate the partition algorithm for finding a labelling partition for a substitution tiling and then contruct the corresponding substitution rule. 
	
	We consider the tiling with interval projection structure introduced in Figures \ref{2to1_interval_staircase}, \ref{2to1_window_eg} and \ref{2to1_applying_g} above. In this example: $$ M= 
	\begin{pmatrix}
		2&1\\1&1 
	\end{pmatrix}
	$$ with $V$ and $W$ the expanding and contracting eigenspaces of $M$. We have, $\Pi_V((1,0)^T) = \tau/\eta$, $\Pi_V((0,1)^T) = 1/\eta$, $\Pi_Y((1,0)^T) = -1$ and $\Pi_Y((0,1)^T) = \tau$. Recall that $\tau$ is the golden ratio and define $\eta = \sqrt{2+\tau}$. The window has length $(2+2 \tau)$ in $Y$ and the staircase has steps $a = (1,-1)^T$, $b = (-1,2)^T$ and $a+b = (0,1)^T$. From Figure \ref{2to1_window_eg} one observes that the window is divided into three intervals with boundary points $\Pi_Y(t+l-b) = \Pi_Y(t + (-1,0)^T ) = \Pi_Y(t) + 1$ and $\Pi_Y(t - a) = \Pi_Y(t + (-1,1)^T) = \Pi_Y(t) + (1-\tau)$.
	
	Our aim is to describe the construction of the substitution rule associated to the subwindow $\ol{\Omega}_{l,t,s} = [\Pi_Y(t) + 1, \Pi_Y(t) + 3) = \Pi_Y(M \Omega_{l,0} + t + s)$ with $s = (-1,0)^T$ of the window $\Omega_{l,t} = [t, t+(2+2 \tau))$, with $l = (-2,2)^T$. See also Figure \ref{2to1_applying_g}. As the substitution rule is independent of $t$, we set $t = (0,0)$ without loss of generality.

	\subsubsection{The partition algorithm} \label{sect:The_partition_eg}
	
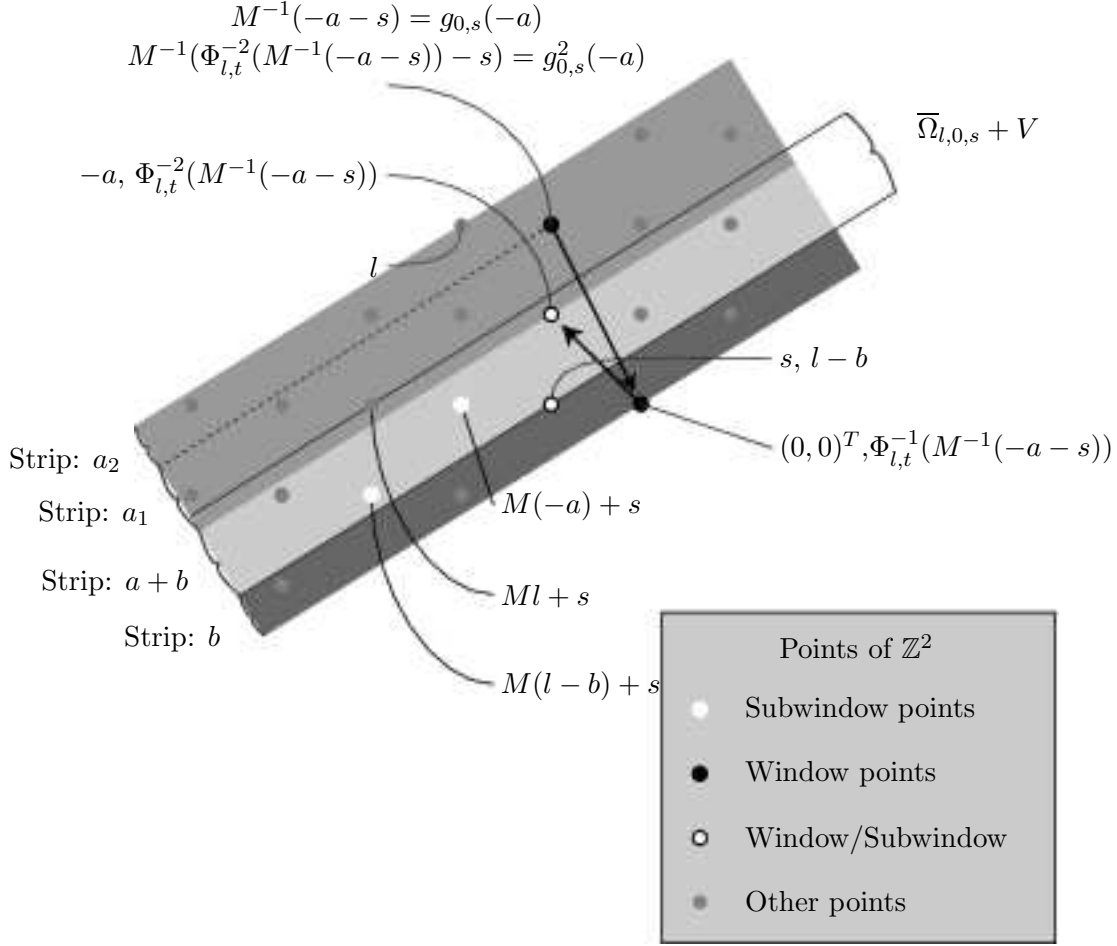
\begin{figure}
		[htp] $$ 
		\begin{xy}
			\xyMarkedImport{} 
			\xyMarkedPos{1}*+!L\txt\labeltextstyle{$\ol{\Omega}_{l,0,s} + V$} 
			\xyMarkedPos{2}*+!R(.8)\txt\labeltextstyle{Strip: $b$} 
			\xyMarkedPos{3}*+!R(.8)\txt\labeltextstyle{Strip: $a+b$} 
			\xyMarkedPos{4}*+!R(.8)\txt\labeltextstyle{Strip: $a_1$} 
			\xyMarkedPos{5}*+!R(.8)\txt\labeltextstyle{Strip: $a_2$} 
			\xyMarkedPos{6}*+!D\txt\labeltextstyle{$M^{-1}(-a-s) = g_{0,s}(-a)$\\ $M^{-1}( \Phi_{l,t}^{-2}(M^{-1}(-a-s)) - s) = g_{0,s}^2(-a)$} 
			\xyMarkedPos{7}*+!R\txt\labeltextstyle{$-a$, $\Phi_{l,t}^{-2}(M^{-1}(-a-s))$} 
			\xyMarkedPos{8}*+!R\txt\labeltextstyle{$l$} 
			\xyMarkedPos{9}*+!L\txt\labeltextstyle{$s$, $l-b$} 
			\xyMarkedPos{10}*+!L\txt\labeltextstyle{$(0,0)^T$,$\Phi_{l,t}^{-1}(M^{-1}(-a-s))$} 
			\xyMarkedPos{13}*+!L\txt\labeltextstyle{$M(l-b)+s$} 
			\xyMarkedPos{12}*+!L\txt\labeltextstyle{$M l+s$} 
			\xyMarkedPos{11}*+!L\txt\labeltextstyle{$M(-a)+s$} 
			\xyMarkedPos{14}*+!L\txt\labeltextstyle{Subwindow points} 
			\xyMarkedPos{15}*+!L\txt\labeltextstyle{Window points} 
			\xyMarkedPos{16}*+!L\txt\labeltextstyle{Window/Subwindow} 
			\xyMarkedPos{17}*+!L\txt\labeltextstyle{Other points} 
			\xyMarkedPos{18}*\txt\labeltextstyle{Points of $\integer^2$} 
		\end{xy}
		$$ \caption[Finding a $2$ to $1$ substitution rule (1)]{Illustration of the construction discussed in Section~\ref{sect:The_partition_eg}. The arrows indicate the action of $\Phi_{l,t}^{-1}$.} \label{2to1_partition_eg} 
	\end{figure}

	The algorithm is illustrated in Figure \ref{2to1_partition_eg}.
	
	\WARMprocessEPS{2to1_partition_eg}{eps}{bb} 
	
	We start with the points $-a$ and $l - b$, which project to $W$ as the boundary points of the partition for tile type, $-a = (-1,1)^T$ and $l - b = (-1,0)^T$. We calculate $g_{l,t,s}(-a)$: the point $\Pi_Y(-a)$ is in $\ol{\Omega}_{l,t,s}$ so $g_{l,t,s}(-a) = M^{-1}(-a-s) = (-1,2)^T$. Next, we calculate $g_{l,t,s}^2(-a)$. We have $\Pi_Y((-1,2)^T) \notin \ol{\Omega}_{l,t,s}$. Moreover $\Phi_{l,t}^{-1}((-1,2)^T) = (0,0)^T$ and $\Pi_Y((0,0)^T) \notin \ol{\Omega}_{l,t,s}$, but $\Phi_{l,t}^{-2}((-1,2)^T) = (-1,1)^T = -a$ projects to $\ol{\Omega}_{l,t,s}$. Hence $-a = g_{l,t,s}^2(-a)$ and the $g_{l,t,s}$-orbit of $-a$ has period two.
	
	We repeat the procedure for $l - b$. We have $\Pi_Y(l -b) \in \ol{\Omega}_{l,t,s}$. Thus $g_{l,t,s}^2(l-b) = M^{-1}(l-b-s) = (0,0)^T$. Subsequently $g_{l,t,s}^2(l-b)(-1,2)^T = g_{l,t,s}(-a)$, considered above. Thus the tail of the $g_{l,t,s}$-orbit of $l-b$ coincides with the $g_{l,t,s}$-orbit of $-a$.
	
	The union of the $g_{l,t,s}$-orbits discussed above contains four points whose $\Pi_Y$-projection yields $G_{l,t,s}$. In Figure \ref{2to1_partition_eg} these points are labelled ``Window points''. By applying $M$ and then translating by $s$, the set of window points are mapped to the set of points labelled ``Subwindow points'' in the figure. The relative positions of the $\Pi_Y$-projection of the subwindow points in $\ol{\Omega}_{l,t,s}$ is in correspondence with the relative positions of the points $G_{l,t,s}$ in $\Omega_{l,t}$.
	
	The points $G_{l,t,s}$ induces a partition of $\Omega_{l,t}$ into four labelling intervals. One interval is associated to the tiles of length $\Pi_V(b)$, one to the tiles of length $\Pi_V(a+b)$, and two to the tiles of length $\Pi_V(a)$. We refer to the corresponding labelled prototiles as $a_1$, $a_2$, $a+b$ and $b$.

	\subsubsection{The replacement rule}
	
\begin{figure}
		[htp] $$ 
		\begin{xy}
			\xyMarkedImport{} 
			\xyMarkedPos{1}*\txt\labeltextstyle{$(a)$} 
			\xyMarkedPos{2}*\txt\labeltextstyle{$(b)$} 
			\xyMarkedPos{3}*\txt\labeltextstyle{$(c)$} 
			\xyMarkedPos{4}*\txt\labeltextstyle{$(d)$} 
			\xyMarkedPos{6}*!R(.9)\txt\labeltextstyle{$\Pi_Y$ of this line\\is the window\\for $a+b$ steps in\\the subwindow} 
		\end{xy}
		$$ \caption[Construction of a $2$ to $1$ substitution rule (2)]{Finding the replacement rule for $\ol{a+b}$. The interpretation of strips and points is as in Figure~\ref{2to1_partition_eg}. The figures illustrate the start of the $\Phi_{l,t}$-orbit of the subwindow of $\ol{\Omega}_{l,t,s}$ associated to the tile $\ol{a+b}$. We note that every iterate of this interval lies entirely within the strip assocated to one of the labelled tiles ($a_1$, $a_2$, $a+b$ or $b$), until the third iterate, in $(d)$, returns to $\ol{\Omega}_{l,t,s}$. The $\Phi_{l,t}$-orbit visits the strips corresponding to the tiles $(a+b)a_1ba_2$ in order. It follows that $\ol{a+b}$ is replaced by $(a+b)a_1ba_2$.} \label{2to1_substitution_eg} 
	\end{figure}
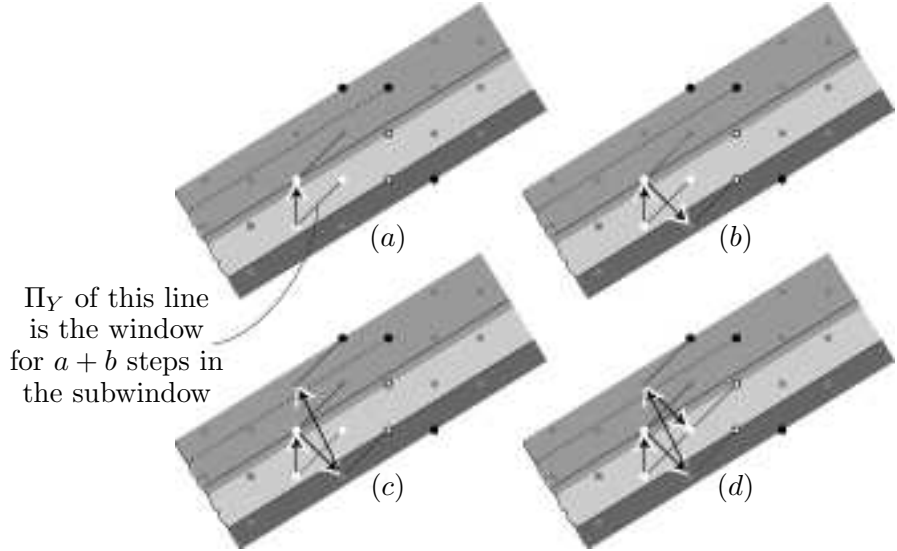

	By following the translations of the labelling intervals of $\ol{\Omega}_{l,t,s}$ induced by $\Phi_{l,t}$, we obtain the replacement rule for the tiling and subwindow considered above. To illustrate this we consider the replacement rule for the tile $\ol{a+b}$ in Figure \ref{2to1_substitution_eg}, showing that $\ol{a+b}$ is replaced by $(a+b)a_1ba_2$. For details we refer to the figure.
	
	\WARMprocessEPS{2to1_substitution_eg}{eps}{bb} 
	
	In a similar fashion one can obtain the replacement rules for the other tiles, yielding the substitution rule: $(a_1\to (a+b)a_1,a_2\to a_1 b a_2, b\to (a+b), a+b \to (a+b)a_1 b a_2)$.

	\section{Substitution rules} \label{sec:proof_subst}

In this section we prove our main results about the characterisation of 
interval projection substitution tilings and their substitution rules. 

The following lemma describes when the set $G_{l,t,s}$ (that has been defined 
in Section~\ref{partition_alg}) is finite. Recall from Lemma~\ref{2to1_finiteifsubst} that this condition is central to the characterisation of substitution rules. 
\begin{lemma}
\label{2to1_finiteifcond}
$G_{l,t,s}$ is finite if and only if $\Pi_Y(l), \Pi_Y(s) \in \rational(\lambda)$.  
\end{lemma}
\begin{proof}
We first show that $G_{l,t,s}$ is finite if $\Pi_Y(l), \Pi_Y(s) \in \Pi_Y(\rational^2)$. Note that $\Pi_Y(\rational^2)=\rational[\lambda]$.  Since only the projections to $Y$ of $l$ and $s$ are meaningful parameters for the tiling, we may assume without loss of generality that $l,s \in \rational^2$.  We show that the set $\{g_{0,s}^n(\Pi_Y(t+l - a))| n \in \nat\}$ is finite. By a similar argument it can be shown that the set $\{g_{0,s}^n(\Pi_Y(t - b))| n \in \nat\}$ is 
finite as well.  Recall that the set $G_{l,t,s} = G_{l,0,s} + t$, so from 
here we only consider $G_{l,0,s}$, without loss of generality.  
 
As each staircase step has positive length and $a + b$ is the largest 
possible step, we have
$$
\Pi_V(z-a+b) \le \Pi_V(\Phi_{l,t}^{-1}(z)) \le \Pi_V(z).
$$
It thus follows that
$$
\Pi_V(z - r (a+b)) \le \Pi_V(\Phi_{l,t}^{-r}(z)) \le \Pi_V(z)~~(r \in \nat).
$$
Now consider the value of $n$ in the definition of $g_{l,t,s}$.  As the steps in the expansion tiling have a maximum length of $\Pi_V(M(a+b))$ and the steps in the original tiling have a minimum length of $\min\{\Pi_V(a), \Pi_V(b)\}$, $n$ has an upper bound $$N = Floor(\frac{\Pi_V(M(a+b))}{\min\{\Pi_V(a), \Pi_V(b)\}}).$$

Now let $\gamma_i = g_{l,0,s}^i(\Pi_Y(t+l - a))$.  From the observations above, and the definition of $\Phi_{l,0}$ and $g_{l,0,s}$ we obtain
$$
\Pi_V(\gamma_{i-1} - N(a+b)) \le \Pi_V(M \gamma_i + s) \le \Pi_V(\gamma_{i-1}),
$$
since $\Pi_V(M \gamma_i + s) = \Pi_V(M g_{l,0,s}(\gamma_i) + s ) = \Pi_V(\Phi_{l,0}^{-n}(\gamma_{i-1}))$ and $n \le N$.  Consequently, as $M$ commutes with 
the projection $\Pi_V$, we have
$$
\Pi_V(M^{-1}(\gamma_{i-1} - N(a+b) - s)) \le \Pi_V(\gamma_i) \le \Pi_V(M^{-1}(\gamma_{i-1} - s)).
$$
By induction on the inequalities above we obtain the following inequality
$$
\Pi_V(\gamma_i) \ge \Pi_V(M^{-i}(\gamma_{0}) - (\sum_{k=0}^i M^{-k}) (N(a + b) + s)).
$$
The matrix $M^{-1}$ acts on $V$ as multiplication by $0 < \lambda^{-1} < 1$, 
so that for all $i$
$$
\Pi_V(\gamma_0 - \frac{1}{1-\lambda} (N(a+b) + s)) \le \Pi_V(\gamma_i) \le \Pi_V(\gamma_0).
$$
Now recall that $s \in \rational^2$, thus $s \in \integer^2/m_s$ for some $m_s \in \nat$.  We observe that for any $n \in \nat$
$$g_{l,0,s}((\integer^2/m_s n)\cap (V + \Omega_{l,0})) \subset ((\integer^2/m_s n)\cap (V + \Omega_{l,0})).$$  
As $\gamma_0 = 0 \in \integer^2$ we thus have 
$$
\gamma_i \in \integer^2/m_s \cap (\Omega_0 + [\Pi_V(\gamma_0 - \frac{1}{1-\lambda} (N(a+b) + s)), \Pi_V(\gamma_0)])
$$
for all $i$.  Thus $\gamma_i$ is contained within a finite subset of the 
lattice $\integer^2/m_s$.  As $l \in \rational^2$, we have 
$l \in \integer^2/m_l$ for some $m_l \in \nat$. 

Using an analogous argument one can show that the points generated by applying $g_{l,0,s}$ to $l$ are also finite.  The set $G_{l,0,s}$ is therefore the projection of the union of two finite sets and is therefore finite.
  
It remains to be shown that whenever $G_{l,0,s}$ is finite we must have $\Pi_Y(l),\Pi_Y(s)\in\rational[\lambda]$.

We have $g_{l,0,s}(l)\in M^{-1} (l-s) + \integer^2$.  As $W$ is an eigenspace of $M$, we have
$$
\Pi_Y(g_{l,0,s}^i(l)) \in \lambda^i \Pi_Y(l-s)+ \Pi_Y(\integer^2).
$$
If the set $G_{l,0,s}$ is finite 
there must exist $i,j \in \mathbb{N}$ with $i\neq j$, such that
$$
(\lambda^i - \lambda^j) \Pi_Y(l-s)\in \Pi_Y(\integer^2).
$$
However, from the properties of $M$ we have $\lambda \in \Pi_Y(\integer^2)$ and it follows that
if $G_{l,0,s}$ is finite we must have $\Pi_Y(l-s)\in\rational[\lambda]$.
The same argument can be repeated, starting with the action of  $g_{l,0,s}$ on $0$ 
(the other extremum of the window), yielding $\Pi_Y(s)\in\rational[\lambda]$. In turn it follows that also
$\Pi_Y(l)\in\rational[\lambda]$
\end{proof}
Our main result is the following Theorem, from which Theorem~\ref{thm:char} 
and Theorem~\ref{thm:subchar} are obtained as immediate corollaries.
\begin{theorem}
\label{thm:main_char}
Let $V$ and $W$ be the tiling and window spaces for the set of interval 
projection tilings $\{\tile_{l,t}|t \in \real^2\}$, with fixed window vector 
$l$ (and window length $\Pi_Y(l)$).  
These interval projection tilings are interval substitution tilings if and 
only if $V$ and $W$ are the eigenspaces of a primitive 
matrix $M\in Gl(2,\integer)$, and $\Pi_Y(l) \in \rational[\lambda]$. 

Each substitution rule is associated to a 
subwindow $\ol{\Omega}_{l,t,s}=\Pi_Y M \Omega_{l,0} + \Pi_Y(s + t)$ with 
$s \in \rational^2$ such that $\ol{\Omega}_{l,t,s} \subset \Omega_{l,t}$.
Substitution rules associated to different choices of $s$ are not patch 
equivalent.
\end{theorem}
\begin{proof}
By Lemma \ref{2to1_finiteifsubst} $G_{l,t,s}$ is finite if and only if there exists a 
substitution rule, and by Lemma \ref{2to1_finiteifcond} the set $G_{l,t,s}$ is finite if and 
only if the conditions set out in the statement of the theorem are satisfied.
\end{proof}

To describe local substitution rules we prove a theorem that is equivalent to Theorem~\ref{thm:local}.  

\begin{theorem}
	\label{thm:main_local}
Consider a local isomorphism class of substitution tilings given by $M$ and $l$ as described in Theorem \ref{thm:main_char}, with substitution rules given by $s$.  
\begin{itemize}
\item If $l \in \integer^2$ the local substitution rules have $s \in \integer^2$
\item If $l - Ml \in \integer^2$ the local substitution rules have $s \in \integer^2$
\item If $l + Ml \in \integer^2$ the local substitution rules have $s \in \integer^2 + l$
\item For all other values of $l$, there are no local substitution rules.
\end{itemize}
\end{theorem}
\begin{proof}
A substitution rule is local if the label of each tile can be 
deduced from the geometric shape of some bounded patch of tiles 
around it. 

Section \ref{partition_defn} shows that the geometric 
shape of a neighbourhood of the tile can be deduced from the projection 
of its vertices to the window. The window admits a hierarchy of partitions, with 
the interval in the partition to which a vertex projects giving the shape of the 
neighbourhood of the tiling near the vertex.  Finer partitions give more information
about the neighbour, in particular giving larger patches of tiling around the vertex.
The boundary points of the union of all these partitions 
form the set $\Pi_Y(\integer^2+t) \cup \Pi_Y(\integer^2+l+t)$. 

A substitution rule is local if the labelling partition is such that its
boundary points precisely lie in the set of boundary point for the partition
that corresponds to the geometric shape of neighbourhoods 
$\Pi_Y(\integer^2+t)\cup \Pi_Y(\integer^2+l+t)$. The boundary points are made up of the projection of the $g_{l,0,s}$ orbits of the points $-a$ and $l - b$.  

Consider the action of $g_{l,0,s}$ on the set $\integer^2 + k$, for some point $k \in \real^2$.  The first part of the function that takes a point to the strip defined by the subwindow adds a vector in $\integer^2$. The second part of the function applies a linear map to expand the subwindow to the full window.  This gives points in the set $$M^{-1}(\integer^2 + k - s) = \integer^2 + M^{-1}(k - s).$$  Thus $g_{l,0,s}(\integer^2 + k) \subset \integer^2 + M^{-1}(k - s)$

Now $-a \in \integer^2$, so $g_{l,0,s}(-a) \in \integer^2 - M^{-1}s$.  This implies that $-M^{-1}s \in \integer^2 \cup (\integer^2 + l)$.  This gives two cases $s \in \integer^2$ and $s \in \integer^2 - M l$, which we will consider separately.

If $s \in \integer^2$ then $$g_{l,0,s}(\integer^2) \subset \integer^2,$$ thus all elements of the orbit of $-a$ correspond to local labels.  Now consider $l - b$.  We have $g_{l,0,s}(l-b) \in \integer^2 + M^{-1}(l - s)$.  As $s \in \integer^2$, there are two possibilities $l \in \integer^2$ and $M^{-1} l \in \integer^2 + l$, the second is equivalent to $l - M l \in \integer^2$.  In both cases, the set $\integer^2 \cup (\integer^2 + l)$ is closed under $g_{l,0,s}$, so the $g_{l,0,s}$-orbits of $-a$ and $l-b$ will only give points in $\integer^2 \cup (\integer^2 + l)$ and the substitution rule is therefore local.  If neither of these conditions are satisfied, then the partition will contain points which are not in $\integer^2 \cup (\integer^2 + l)$ and the substitution rule will not be local.

Now consider $s \in \integer^2 - M l$, where $g_{l,0,s}(-a) \in \integer^2 + l$.  We must therefore consider the points in $g_{l,0,s}(\integer^2 + l)$, to obtain further points on the orbit of $-a$ and the orbit of $l-b$.  We have $$g_{l,0,s}(\integer^2 + l) \subset \integer^2 + M^{-1}(l-s) = \integer^2 + M^{-1} l + l.$$  Again this leads us to consider the two cases $l \in \integer^2$ and $M^{-1} l + l \in \integer^2$.  By an analogous argument to the one given above, the substitution rule is local if and only if one of these conditions holds.  

We reorder the conditions obtained above to emphasize the length of the window (and thus the local isomorphism class of the tiling) first, yielding the following cases:
\begin{enumerate}
\item $l \in \integer^2$ and $s \in \integer^2 \cup \integer^2 - M l$
\label{enum:localcases_proof_1}
\item $l - M l \in \integer^2$ and $s \in \integer^2$
\label{enum:localcases_proof_2}
\item $l + M l \in \integer^2$ and $s \in \integer^2 - M l$
\label{enum:localcases_proof_3}
\end{enumerate}

These may be further simplified to yield the conditions given in the statement of the theorem.  In the first case the two sets that can contain $s$ are the same since $l \in \integer^2$.  In the third case we may add $l + M l \in \integer^2$ to the set for $s$, to obtain $s \in \integer^2 + l$
\end{proof}

Finally we show how this result is equivalent to Theorem \ref{thm:local}.

\begin{proof}[Proof of Theorem \ref{thm:local}]
By Theorem \ref{thm:main_local} we have three necessary and sufficient conditions for the existence of local substitution rules:
\begin{enumerate}
\item $l \in \integer^2$ and $s \in \integer^2$
\label{enum:localcases2_proof_1}
\item $l - M l \in \integer^2$ and $s \in \integer^2$
\label{enum:localcases2_proof_2}
\item $l + M l \in \integer^2$ and $s \in \integer^2 + l$
\label{enum:localcases2_proof_3}
\end{enumerate}
The first case holds precisely if both the second and third case hold, so we need only consider those.  

If $\det M=1$, then $d_- = |\Pi_Y(s)|$ and $d_+ = |\Pi_Y(l - M l - s)|$. If $\det M=-1$ then $d_- = |\Pi_Y(s + M l)|$ and $d_+ = |\Pi_Y(l - s)|$.  There thus remain four cases to be considered, 
the second and third cases from above with $\det M$ either positive or negative.
Firstly if $\det M=1$, the second case we have $(1-\det(M)/\lambda)|\Omega|, d_+, d_- \in \Pi_Y(\integer^2)$, as $\Pi_Y(l) = |\Omega|$ and $\Pi_Y(M l) = \det(M)/\lambda |\Omega|$.  In the third case we have $(1+\det(M)/\lambda)|\Omega|\in \Pi_Y(\integer^2)$, $d_- \in \Pi_Y(\integer^2 + l)$ and $$d_+ \in \Pi_Y(l - M l - \integer^2 - l) = \Pi_Y(\integer^2 + l),$$
since $l + M l \in \integer^2$.

Finally consider the case that $\det M=-1$. In this case $\Pi_Y(M l)$ has the opposite orientation of $\Pi_Y(l)$.  Hence, in the second case it follows that $|\Omega| + |M \Omega| \in \Pi_Y(\integer^2)$, $d_- \in \Pi_Y(\integer^2 + Ml) = \Pi_Y(\integer^2 + l)$ and $d_- \in \Pi_Y(\integer^2 + l)$.  In the third case we have $|\Omega| - |M \Omega| \in \Pi_Y(\integer^2)$, $d_- \in \Pi_Y(\integer^2 + l + M l) = \Pi_Y(\integer^2)$ and $d_+ \in \Pi_Y(\integer^2 + l - l) = \Pi_Y(\integer^2)$. 

We may therefore summarize the cases as in Theorem \ref{thm:local}:
\begin{itemize}
\item $(1-\det(M)/\lambda)|\Omega|, d_+,d_-\in \Pi_Y(\integer^2) = \integer[\lambda]$, 
\item $(1+\det(M)/\lambda)|\Omega|, d_+ - |\Omega|,d_- -|\Omega|\in \integer[\lambda]$.
\end{itemize}
\end{proof}
	
\section{Tilings fixed by a substitution rule}
\label{sec:fixed}

In this section we address the following question: given a substitution rule, which
tilings are fixed, ie mapped to (a translated copy of) themselves, by this substitution? 
Recall that in general the predecessor of a substitution tiling is only locally isomorphic to 
this tiling, and often not a (translated) copy.

Because of corresponding interests in the literature~\cite{Crisp:SICS, Komatsu:SIBS, Yasutomi:OSSWA}, we consider not only tilings of the 
entire line, but also tilings of the positive half-line $\real^+\cup\{0\}$ (with boundary vertex at $0$).

The results are summarised in the following theorem, which contains the result announced in 
Theorem~\ref{introthm:fix_tilings}.
\begin{theorem}
\label{thm:fixed}
Consider a substitution rule $\sigma$ with interval projection structure. Let 
$s,l\in\rational^2$ denote the associated subwindow translation and window vectors. 
Let $\tilde{\lambda}$ ($=\det(M)/\lambda$) denote the contracting eigenvalue of the
primitive matrix $M\in Gl(2,\integer)$, associated to the substitution rule. Then 

the number of tilings of the line and half-line that are fixed by the substitution rule are as listed in 
Table~\ref{table_two_sided} and Table~\ref{table_one_sided}. These numbers are functions of
$s$, $l$, $M$, and the invariant window interval (IWI) 
$$[\frac{\Pi_Y(s)}{\tilde\lambda -1},\frac{\Pi_Y(s)}{\tilde\lambda -1}+\Pi_Y(l)] \subset W.$$
\end{theorem}

\begin{table}[htp]

\centering

\begin{tabular}{|c||c|c|} \hline

 & $\det(M)=1$ & $\det(M)=-1$\\ 

\hline \hline 
$s\not\in \integer^2\cup \integer^2-(M-I)l$
& $\det(M-I)$ & $\det(M-I)$\\

\hline
$s\in \integer^2\cup \integer^2-(M-I)l$
& $\det(M-I) + 1$ & $\det(M-I) - 1$\\ 

\hline

\end{tabular}
\caption{The number of tilings of the line that are fixed by a substitution rule with an interval 
projection structure.  $M\in Gl(2,\integer)$ is the primitive matrix and  $s,l\in\rational^2$ 
are the subwindow translation vector and window vector associated with the substitution rule.  
$I$ denotes the $2\times 2$ identity matrix. }
\label{table_two_sided}
\end{table}

\begin{table}[htp]

\centering

\begin{tabular}{|c||c|c|} \hline

 & $\det(M)=1$ & $\det(M)=-1$\\ 

\hline \hline 

IWI non-singular & $1$ & $1$\\

\hline

IWI singular & $1$ & -\\

$0$ on boundary & & \\
\hline

IWI singular & $2$ & $0$\\

$0$ in interior & & \\

\hline

\end{tabular}

\caption{The number of tilings of the half-line $\real^+\cup\{0\}$
that are fixed by a substitution rule with an interval 
projection structure with associated primitive matrix $M\in Gl(2,\integer)$, subwindow translation 
vector and window vector $s,l\in\rational[\lambda]$, and invariant window interval (IWI)
$[\frac{\Pi_Y(s)}{\tilde\lambda -1},\frac{\Pi_Y(s)}{\tilde\lambda -1}+l] \subset W$.
The IWI is singular if and only if $s\in (M-I)\integer^2\cup (M-I)(\integer^2-l)$.
}

\label{table_one_sided}
\end{table}

Key to the proof of Theorem~\ref{thm:fixed} is the observation that the substitution rule induces a simple map on the set of windows taking the the window of the original tiling $\mathcal{T}$ to the window of the tiling $\sigma{\mathcal{T}}$ obtained by
application of the substitution rule. 

\begin{lemma}
\label{lem_subst_image}
Consider a tiling $\tile_{l,t}$ with an interval projection structure and
window base vector $t\in \real^2$, that admits a substitution rule $\sigma$ as in Theorem~\ref{thm:fixed}.
Then  $\sigma(\tile_{l,t}) = \tile_{l,M t - s}$.
\end{lemma}

\begin{proof}
We recall that the subwindow of $\Omega_{l,r}$ associated with the 
substitution rule is $\ol{\Omega}_{l,r,s}:=\tilde\lambda\Omega_{l,0}+r + s$.
Then $M^{-1}\ol{\Omega}_{l,r,s}=\Omega_{l,M^{-1}(r+s)}=\Omega_{l,t}$
is the window associated to the predecessor. As $r=M t- s$ the substitution 
rule thus induces the transformation $t\to M t- s$. 
\end{proof}

Finally we prove Theorem~\ref{thm:fixed}.

\begin{proof}[Proof of Theorem \ref{thm:fixed}] 
We first consider tilings of the entire line $\real$.  We note that two projection tilings are the same up to translation if the projections of the associated window base vectors to $Y$ differ by a vector in $\Pi_Y(\integer^2)$.  

Using the result of Lemma~\ref{lem_subst_image} we thus have
\begin{equation}
\Pi_Y(t+z)=\Pi_Y( M t - s)\label{transproj}
\end{equation}
for some $z \in \integer^2$.  

As before, without loss we consider $s\in\rational^2$.
Recall that if $t\in\real^2$ is a window base vector for a projection tiling, then 
any vector in $\{t+v~|~v\in V\}$ is also a window base vector for the same tiling.
Since $M\in Gl(2,\integer)$, we may in fact without loss of generality restrict the solution
set of (\ref{transproj}) to $t\in\rational^2$ and consequently replace this equation by

\begin{equation}\label{trans2}
t = (M - I)^{-1}(s + z),~~z\in \integer^2
\end{equation}
where $I$ denotes the $2\times2$ identity matrix.  

Viewing (\ref{trans2}) to define $t$ as a function of $z$, (equivalence classes of) 
tilings that are invariant under this substitution rule are represented by the elements of 
$\{t(z)\mod\integer^2~|~z\in\integer^2\}$. 

The number of different (closed) window intervals that are invariant under a substitution rule 
is thus equal to
$$
|\{{t}(z)\mod\integer^2~|~z\in\integer^2\}|=|\left((M - I)^{-1}\integer^2\right)/ \integer^2|=|\det(M - I)|.
$$
Note that $\det(M-I)=1-T+D$, where $T={\rm tr}M$, $D=\det M$ and $\lambda^2-T\lambda+D=0$.

In the argument above we find the number of closed window intervals that are invariant under the 
induced action of a substitution rule. In order to make the relation with tilings we need to 
recall that the window is a half-open interval, which provides the window not only with
a position and length, but also with an orientation. 

We say that a tiling with interval projection structure (and its window) is \emph{singular} if $\partial\Omega_{l,t}+V\cap \integer^2\neq \emptyset$ (or equivalently if $\partial\Omega_{l,t}\cap \Pi_Y(\integer^2)\neq\emptyset$) and non-singular otherwise. 

The base vector $t$ of a singular window satisfies
\begin{equation}
\Pi_Y(t)\in \Pi_Y(\integer^2)\cup\Pi_Y(\integer^2)-\Pi_Y(l).
\label{sing}
\end{equation}
Since $\Pi_Y(l)\in\rational[\lambda]$ we can without loss choose $l\in\rational^2$ and subsequently 
solve for $t\in\rational^2$ in the case of singular tilings, rewriting (\ref{sing}) as
$$
t\in \integer^2\cup \integer^2-l.
$$
It thus follows that within a local isomorphism class there are one ($l\in\integer^2$) or 
two ($l\not\in\integer^2$) different tilings (up to translation) that are singular. 

Consequently, there exists a singular tiling that is candidate to be invariant 
(a $t\in\integer^2\cup\integer^2-l$ that satisfies (\ref{trans2})) if and only if
$$
s\in \integer^2\cup \integer^2-(M-I)l.
$$

If $s\not\in \integer^2\cup \integer^2-(M-I)l$ there exists
no singular invariant tiling and
the orientation of the window needs no attention as the boundary of the window of an invariant tiling
cannot intersect the lattice $\integer^2$. For such substitution rules the number of different invariant 
tilings is thus exactly equal to $|\det(M-I)|$. 

In the singular case it is important to compare the orientation of the window before and after 
substitution. Note that for each $s\in \integer^2\cup \integer^2-(M-I)l$ there is precisely one choice 
of $t$ (mod $\integer^2$) satisfying (\ref{sing}) and thus precisely one singular tiling that is a candidate to be invariant under the substitution rule.

If $\det M=1$ the windows before and after the geometric transformation associated with
substitution have the same orientation. For each singular 
window base vector $t$, we thus find two different singular invariant tilings: one for each choice of the 
side of closure of the window. Consequently, in this case the number of invariant tilings is 
$|\det(M-I)|+1$. When $\det M=-1$ the orientation of the singular tiling changes so the 
substitution rule cannot fix it and the total number of invariant tilings is $|\det(M-I)|-1$. 

The results are summarized in Table~\ref{table_two_sided}.

It now remains to analyse the number of fixed tilings of the half-line $\real^+\cup\{0\}$.
In this case there is no way to allow for translations, so that invariance implies
$$
t = M t - s ~\Leftrightarrow~ t = (M - I)^{-1}s,
$$
with $t,s\in\rational^2$. The corresponding invariant window interval (IWI) is 
$[\frac{\Pi_Y(s)}{\tilde\lambda -1},\frac{\Pi_Y(s)}{\tilde\lambda -1}+l] \subset W$.

We first note that $0\in W$ lies inside the IWI.  In fact, $0$ is the unique vertex that is also a 
vertex for all the expansion predecessors.  Hence $0$ must be the boundary vertex.  

If the window translation vector $s$ of a substitution rule satisfies
$s\not\in (M-I)\integer^2 \cup (M-I)(\integer^2 - l)$ then the IWI is non-singular and there 
is a unique invariant tiling of the half-line $\real^+\cup\{0\}$.

If the window is singular then we must consider the sign of determinant of the matrix $M$
and the location of (projections) of the singular lattice point(s). If the tiling admits
a vertex that is the projection of a singular lattice point (lying on the boundary of 
$V+\Omega_{l,t}$) then we call this vertex singular.

If the IWI is singular, but the tiling does not contain a singular vertex (as all singular lattice
points project to the negative part of the half-line), there is no interference with
the tiling and thus there is one unique invariant tiling of the half-line $\real^+\cup\{0\}$.

If $0$ is at the boundary of the IWI it is a singular vertex. As $0$ serves as the boundary point 
we have no choice with regard to the side on which to close the window. Moreover it turns out that 
$0$ can be on the boundary only if $\det(M)=+1$. Consequently, there is one unique invariant 
tiling of the half-line $\real^+\cup\{0\}$.

If the tiling for the IWI has a singular vertex that does not lie on the boundary we have to distinguish
two cases, depending on the sign of $\det(M)$. If the sign is positive there are two invariant tilings,
corresponding to the two choices of the side where the window is closed. If the sign is negative there are
no invariant tilings.

The results are summarized in Table~\ref{table_one_sided}.
\end{proof}
	
\section{Interval exchange transformations}\label{sec:iet}
	
Our results on the existence of substitution for one-dimensional tilings with an interval projection structure are intimately related to renormalization properties of two and three-\emph{Interval Exchange Transformations} (IETs). In this section we discuss this relationship and reformulate our main result in terms of IETs.

IETs are defined as follows. Let $I\subset\mathbb{R}$ be a bounded half open interval that is partitioned into a set of half open subintervals 
$\{I_\alpha~|~\alpha\in\mathcal{A}\}$ indexed by some alphabet $\mathcal{A}$ of $d\geq2$ symbols. An interval exchange transformations is a bijection from $I$ to $I$ which is a constant translation on each of the subintervals $I_\alpha$, thus inducing a permutation of the subintervals covering $I$.
An interval exchange transformation $f$ is determined by combinatorial and 
metrical data as follows:
\begin{itemize}
\item[(1)]
A pair $\pi=(\pi_0,\pi_1)$ of bijections 
$\pi_\varepsilon:\mathcal{A}\to\{1,\ldots,d\}$, with $\varepsilon\in\{0,1\}$, 
describing the ordering of the subintervals $I_\alpha$ before and after the transformation is applied. One may represent this permutation as
\[
\pi=\left(\begin{array}{cccc}
\alpha_1^0&\alpha_2^0&\ldots&\alpha_d^0\\
\alpha_1^1&\alpha_2^1&\ldots&\alpha_d^1
\end{array}\right),
\]
where $\alpha_j^\varepsilon=\pi_\varepsilon^{-1}(j)$ and $j\in\{1,2,\ldots,d\}$.
\item[(2)] A vector $\mu=(\mu_\alpha)_{\alpha\in\mathcal{A}}$ with positive entries, where $\mu_\alpha$ denotes the length of the subinterval $I_\alpha$.\end{itemize}
We call $p=\pi_1\circ\pi_0^{-1}$ the \emph{monodromy invariant} of $\pi$. It should be noted that we can normalize the combinatorial data by choosing $\mathcal{A}=\{1,2,\ldots,d\}$ and $\pi_0=\id$ so that $\pi_1=p$ and $\pi=(\id,p)$.
We denote an interval exchange transformation (IET) $f$ by the pair $f=(\pi,\mu)$.

The study of interval exchange transformations has a rich history, mainly motivated by the study of translation surfaces. For more details, we recommend the recent survey by Viana \cite{Viana:ETIEM}.

IETs appear in the study of one-dimensional tilings with an interval projection 
structure, via the map $\Pi_Y\circ \Phi_{l,t}$ whose constituents have been 
introduced in Section~\ref{2to1_IPT}. Orbits of this map represent 
projections of subsequent vertices of the staircase to the window 
$\Omega_{l,t}\subset Y$. From Definition~\ref{defn:Phi} it follows that
$\Pi_Y\circ\Phi_{l,t}$ is an IET on two or three intervals covering the window
$\Omega_{l,t}$.  Such IETs are considered in detail in \cite{Ferenczi:SOTI1, Ferenczi:SOTI2,Ferenczi:SOTI3}.

The type of IETS arising in this context of Section~\ref{2to1_IPT} are precisely
the irreducible two and three-IETs. An IET $f$ on an interval $I$ is called 
irreducible if there does not exist any subinterval $\ol{I}$ of $I$ for which 
$f|_{\ol{I}}$ is an IET on $\ol{I}$. These IETs may be represented by the following 
two-parameter family of maps. Let $I=[0,1)$ denote the half open unit 
interval and $a,b\in I$ two points such that $0<a\leq b<1$. Then we define
\[
f_{a,b}(x)= \left\{ 
	\begin{array}{ll} x + 1 - a & \text{if } x\in [0,a),\\
		x + 1 - a - b & \text{if }x\in [a,b),\\
		x - b & \text{if }x\in [b,1),
	\end{array}
\right. 
\]	
The map $f_{a,b}$ is a two-IET if $a=b$ and a three-IET otherwise.
	
In the study of IETs, considerable attention has been devoted to the return 
map $\hat{R}_{\ol{I}}f$ that an IET $f=(\pi,\mu)$ induces on a half open 
subinterval $\ol{I}\subset I$. It turns out that the return map is again an IET,
with permutation $\ol{\pi}$ of a partition $\ol{\mu}$
of $\ol{I}$ that is implicitly defined by the 
dynamics of $f$. In order to compare the IET $f$ on $I$ with the return map that 
it induces on $\ol{I}\subset I$, it is natural to identify $\ol{I}$ with $I$ 
(by means of a uniform expansion with factor $|I|/|\ol{I}|$, an appropriate 
translation, and - in case the orientations of $I$ and $\ol{I}$ are opposite - a 
reflection). We refer to this mapping between IETs as the 
\emph{Rauzy-Veech renormalization} induced by $\ol{I}$, after \cite{Rauzy:EDTI,Veech:IET}:
\[ R_{\ol{I}}(\pi,\mu)=(\ol{\pi},\ol{\mu}).\]

More information about the dynamics of the original IET is retained by
considering an extension of the Rauzy-Veech 
renormalization that includes suspension data. In addition to considering 
the induced IET on the subinterval $\ol{I}$, one may choose to
keep track of the paths that
suspended orbits take in $I$, while being taken by iteration of $f$ 
from $\ol{I}$ through $I\setminus \ol{I}$ before returning to $\ol{I}$. 
The suspension data may be represented by an element of the monoid 
$\mathcal{A}^*$ consisting of all finite words in the alphabet $\mathcal{A}$ 
that is used to label the subintervals $I_\alpha$ of $I$. 
We thus define the Rauzy-Veech renormalization on the IET including suspension 
data as 
\[
\mathcal{R}_{\ol{I}}(\pi,\mu,v)=(\ol{\pi},\ol{\mu},\ol{v}),
\]
where the suspension data are represented by $v=(v_1,\ldots,v_n)$ 
with $v_i$ an element of the monoid $A^*$.

Recall that a morphism, or substitution rule, 
$\sigma:\mathcal{A}^*\to\mathcal{A}^*$
associates to each element of the alphabet $\mathcal{A}$ a finite word over 
$\mathcal{A}$ and acts on words in  
$w(\alpha_1,\alpha_2,\ldots,\alpha_d)\in\mathcal{A}^*$ 
as 
\[
\sigma w(\alpha_1,\alpha_2,\ldots,\alpha_d)=w(\sigma(\alpha_1),\sigma(\alpha_2),\ldots,\sigma(\alpha_d)).
\]

We address the question for which IET $(\pi,\mu)$ and
subinterval $\ol{I}$ there exists a refinement $\tilde\mu$ of the partition
$\mu$, on which $\pi$ induces a permutation $\tilde\pi$, so that 
the Rauzy-Veech renormalization mapping  fixes 
$(\tilde\pi,\tilde\mu)$ and induces a morphism $\sigma$ on the corresponding
suspension data $\tilde v=(\tilde v_1,\ldots,\tilde v_{\tilde d})$
with $\tilde v_i\in\tilde{\mathcal{A}}^*$ where 
$\tilde{\mathcal{A}}$ is an alphabet with $\tilde{d}$ element:
\begin{equation}\label{RVren}
\mathcal{R}_{\ol{I}}(\tilde\pi,\tilde\mu,\tilde v)=
(\tilde\pi,\tilde\mu,\sigma \tilde v).
\end{equation}

\begin{theorem}\label{thm:iet}
Let $f_{a,b}=(\pi,\mu)$ be a two or three-interval exchange transformation  
on the unit interval $I$ with discontinuity points $0<a\leq b<1$. Then
there exists a refinement $\tilde\mu$ of $\mu$, with permutation $\tilde\pi$ 
trivially induced from $\pi$, such that the Rauzy-Veech renormalization with 
respect to a subinterval $\ol{I}$ induces a morphism $\sigma$ on the suspension data 
$\tilde v$, as in (\ref{RVren}), if and only if there exists a quadratic unit $\lambda$, 
                    such that 
\begin{itemize}
\item $|\lambda|<1$,
\item its algebraic conjugate is positive, 
\item      $\lambda(\integer+\integer (a-1)/b)=\integer+\integer (a-1)/b$,
\end{itemize}        
and $\ol{I}=t+\lambda I \subset I$ with $t\in\rational[\lambda]$.
\end{theorem}
\begin{proof}
This theorem is a consequence of Theorem~\ref{thm:char} and 
Theorem~\ref{thm:subchar}. From the discussion in 
Section~\ref{partition} and Section~\ref{sec:proof_subst} it follows that the
existence of a substitution rule for a one-dimensional tiling with interval
projection structure is equivalent to the existence of an induced morphism
on the suspension date of Rauzy-Veech renormalization for a two or three-IET.
The subinterval $\ol{I}$ corresponds to the subwindow associated to
the predecessor tiling and the refinement $\tilde\mu$ to the labelling 
partition of the window, and the morphism $\sigma$ to the substitution rule.         

Let us consider an IET with interval length $\gamma>0$ and directed
one-dimensional vectors $\alpha <0<\beta$, satisfying 
$|\alpha|+|\beta|\geq\gamma$, as considered in 
Section~\ref{partition}. This corresponds to an IET on the unit
interval $I$ with discontinuity points $0<a\leq b<1$, where 
$a=(\gamma+\alpha)/\gamma$ and $b=\beta/\gamma$.

In the setting of Section~\ref{partition}, $\alpha$ and $\beta$ are 
considered as projections (to $Y$) of generators of a two-dimensional 
lattice isomorphic to $\integer^2$. The projection of this lattice to 
$Y$ is $\alpha \integer + \beta \integer$. 
If there exists a lattice automorphism that acts as a uniform contraction 
by a factor $\lambda$ on $Y$ and as a uniform expansion by a positive
factor containing the 
interval, it follows that $\lambda$ must be a quadratic unit, and the 
point set $\alpha \integer+\beta \integer$ must 
be invariant under multiplication by 
$\lambda$, which implies that
\[
\Pi_Y(\integer^2)(\integer b+\integer(a-1))=\integer b+\integer(a-1).
\]
This implies in particular that there exist integers $x,y$ such that
$\lambda=x+y(a-1)/b$ and hence $\lambda\in \integer+\integer(a-1)/b$, 
but also that $\lambda=xb/(a-1)+y$ so that $(a-1)/b\lambda \in \integer+
\integer(a-1)/b$ which in turn implies (as $\lambda\not\in\integer$)  
that $(a-1)/b$ must be a quadratic algebraic number.

On the window and subwindow $\ol{I}$ we have the following conditions from 
the corresponding conditions for substitution rules in Theorem~\ref{thm:char}
and Theorem~\ref{thm:subchar}. First, the length of the window interval
$\gamma$ must be in $\rational[\lambda]$, and second, so must the 
distances between any of the boundaries of the subwindow interval to 
any boundary of the orginal window interval. In terms of the normalized IET
this implies that the distances between boundary points of $\ol{I}$ to any of 
the boundaries of $I$ must be in $\rational[\lambda]$, i.e.
\[
\partial \ol{I}\subset\rational[\lambda]\cap I.
\]
\end{proof}
Note that the conditions on $\lambda$ are necessary and sufficient conditions for the 
existence of a subinterval $\ol{I}$ on which the Rauzy-Veech renormalization induces
a substitution rule on the suspension data, in correspondence with Theorem~\ref{thm:char}. The final condition corresponds to the one set out in 
Theorem~\ref{thm:subchar} characterizing the countably infinite number of
appropriate choices for the subinterval $\ol{I}$.

Theorem~\ref{thm:iet} relates to some results on renormalization by 
\cite{Boshernitzan:AEOLT}, where Rauzy-Veech renormalization on specific
subintervals is considered, without considering the corresponding action on 
suspension data.

We finally illustrate the Rauzy-Veech renormalization on the Fibonacci IET.

\WARMprocessEPS{Fibb_middle_int_exch}{eps}{bb}
	\begin{figure}[htp]
	\centerline{\begin{xy}
	\xyMarkedImport{}
	\xyMarkedPos{1}*!D\txt\labeltextstyle{$a$}
	\xyMarkedPos{2}*!D\txt\labeltextstyle{$b$}
	\xyMarkedPos{3}*!D\txt\labeltextstyle{$\ol{b}$}
	\xyMarkedPos{4}*!D\txt\labeltextstyle{$\ol{a}$}
	\xyMarkedPos{6}*!D\txt\labeltextstyle{$a_1$}
	\xyMarkedPos{7}*!D\txt\labeltextstyle{$a_2$}
	\xyMarkedPos{8}*!D\txt\labeltextstyle{$b$}
	\xyMarkedPos{9}*!D\txt\labeltextstyle{$\ol{b}$}
	\xyMarkedPos{11}*!D\txt\labeltextstyle{$\ol{a}_2$}
	\xyMarkedPos{12}*!D\txt\labeltextstyle{$\ol{a}_1$}
	\xyMarkedPos{13}*!D\txt\labeltextstyle{$a_1$}
	\xyMarkedPos{14}*!D\txt\labeltextstyle{$a_2$}
	\xyMarkedPos{15}*!D\txt\labeltextstyle{$b_1$}
	\xyMarkedPos{16}*!D\txt\labeltextstyle{$b_2$}
	\xyMarkedPos{17}*!D\txt\labeltextstyle{$\ol{b}_2$}
	\xyMarkedPos{18}*!D\txt\labeltextstyle{$\ol{b}_1$}
	\xyMarkedPos{19}*!D\txt\labeltextstyle{$\ol{a}_2$}
	\xyMarkedPos{20}*!D\txt\labeltextstyle{$\ol{a}_1$}
	\xyMarkedPos{27}*!D\txt\labeltextstyle{(i)}
	\xyMarkedPos{28}*!D\txt\labeltextstyle{(ii)}
	\xyMarkedPos{29}*!D\txt\labeltextstyle{(iii)}
	\end{xy}}
	\caption[Central Fibonacci]{Illustration of the iterative construction of a 
        refined partition of
        the interval $I$ for which the Rauzy-Veech renormalization on the central 
        subinterval $\ol{I}$ of length $\frac{\sqrt{5}-1}{2}=-\lambda$ induces a 
        substitution rule on the suspension data for the Fibonacci IET, defined in 
        Example~\ref{ex:nonlocalfibdet}.  The suspension data is 
        represented as a graph over the partition of the interval. Because 
        the contraction factor $\lambda<0$, in comparing the IET on $I$ to the induced
        IET on $\ol{I}$, we need not only rescale and translate, but also perform a
        reflection in the line. We start with the Fibonacci IET, schematically 
        indicated in the left graph under (i), and the induced IET on $\ol{I}$
        indicated in the right graph. It follows that the Rauzy-Veech renormalization
        fixes the IET $R_{\ol{I}}(\pi,\mu)=(\pi,\mu)$, but it does not induce a 
        morphism on the suspension data, since in the graph over the longest 
        subinterval of $\ol{I}$ one distinguishes two different graphs: $a b$ and $b a$,
        corresonding to different suspensions. In order for the renormalization to
        induce a morphism on the suspension data, corresponding to each marked
        subinterval of $\ol{I}$ there should be unique suspension data. We thus need at
        least make a refinement of the $a$ subinterval of $I$ into two subintervals
        which we label $a_1$ and $a_2$ (left of (ii)). If we examine the tower on 
        the induced partition $\{\ol{I}_{\ol{a}_1},\ol{I}_{\ol{a}_2},\ol{I}_{\ol{b}}\}$ of $\ol{I}$, we observe that
        the graphs over $\ol{I}_{\ol{a}_1}$ and $\ol{I}_{\ol{a}_2}$ are well defined, but that there
        is a problem with the suspension over $\ol{I}_{\ol{b}}$. We subsequently consider
        the additonal refinement depicted at the left of (iii), where the subinterval
        $I_b$ is partitioned in two subintervals $I_{b_1}$ and $I_{b_2}$. With this 
        partition, we observe at the right figure under (iii), that the suspension
        data on the induced IET on $\ol{I}$ now is unique on each labeled subinterval
        of $\ol{I}$. From the graph we read the morphism 
        $
        \sigma (a_1,a_2,b_1,b_2)=(b_1a_1,a_2b_2,a_2,a_1).
        $  
        This is precisely the substition rule mentioned in 
        Example~\ref{ex:nonlocalfibdet}. }
	\label{Fibb_middle_int_exch} 
	\end{figure}
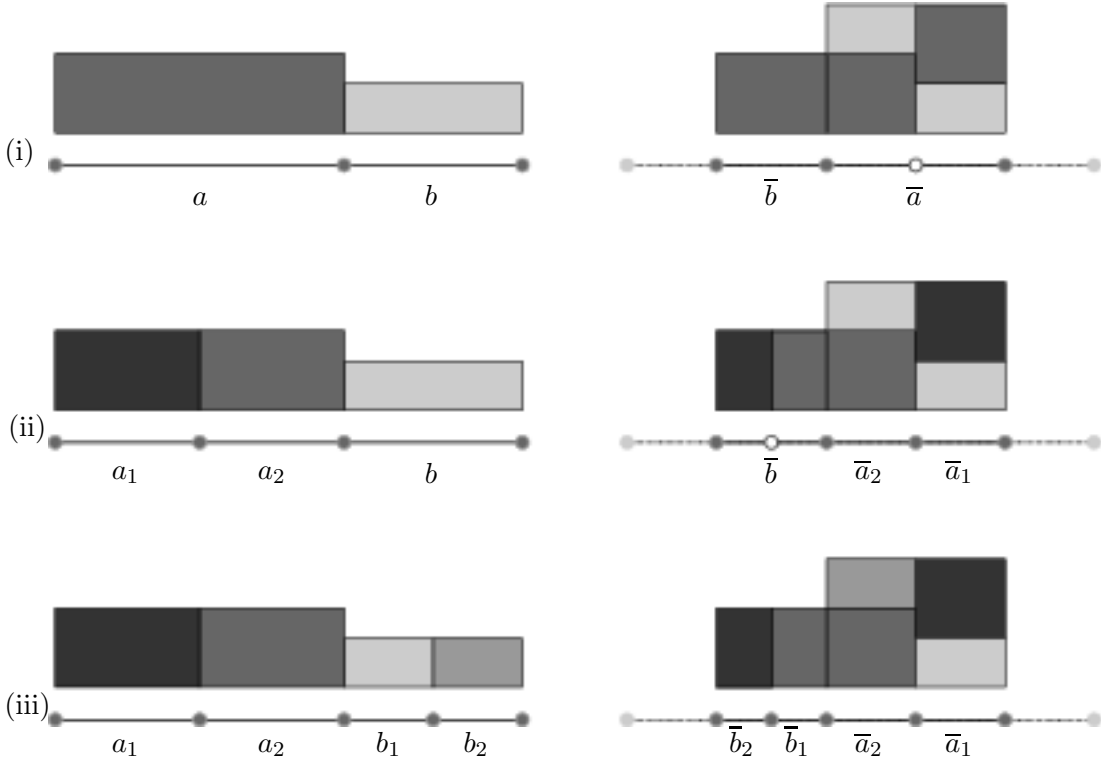
\begin{example}[Rauzy-Veech renormalization of the Fibonacci IET]
\label{ex:nonlocalfibdet} 
To illustrate the connection between renormalization of IETs and substitution 
rules for one-dimensional tilings with an interval projection structure, we 
consider the Fibonacci IET, whose orbits correspond to Fibonacci tilings.	
We focus in particular on the derivation of the non-local substitution rule
for the Fibonacci tilings, presented in Example~\ref{ex:nonlocalfib}. 

The Fibonacci two-IET on $I=[0,1)$ is defined as $f_{a,b}$ with $a=b=
\frac{\sqrt{5}-1}{2}$. Choosing the quadratic unit $\lambda=-a$, the conditions on $\lambda$ of Theorem~\ref{thm:iet} are satisfied.
We choose $\ol{I}=\lambda I+t$ as the central subinterval of the unit interval $I$, 
with $\lambda=-\frac{\sqrt{5}-1}{2}$ a quadratic unit and 
$t=\frac{1+\sqrt{5}}{4}\in\rational[\lambda]\setminus\integer[\lambda]$.
Satisfying hereby the final condition of Theorem~\ref{thm:iet}, we thus expect 
to have some refinement $\tilde\mu$ of the partition $\mu$ for which the 
Rauzy-Veech renormalization to $\ol{I}$ induces a morphism (substitution rule) to 
the suspension data.
 
In Figure~\ref{Fibb_middle_int_exch} we illustrate how the refinement of 
the partition can be constructed by 
examining the return map to $\ol{I}$, including suspension data.
We note that due to the fact that $\lambda<0$, the orientation of 
the half open interval $\ol{I}$ is opposite to that of $I$. This implies that
in the identification of the IET induced on $\ol{I}$ by the return map with $I$, we 
need to take into consideration a reflection. As the boundary points of
$\ol{I}$ lie in $\rational[\lambda]\setminus\integer[\lambda]$, the substitution rule
is non-local, in the sense of the discussion in Section~\ref{sect:Main_results}.
\end{example}
	
\subsection*{Acknowledgements}
EOH and JSWL gratefully acknowledge the support of the UK Engineering and 
Physical Sciences Research Council (EPSRC).

	\nocite{Baake:DIMQ}
	
\bibliography{Quasicrystals}
	\bibliographystyle{amsalpha} 
\end{document}